\documentclass[authoryear, 10pt]{elsarticle}
\usepackage{hyperref}

\usepackage{url}
\usepackage{empheq}
\allowdisplaybreaks

 \DeclareMathOperator{\sign}{sign}

\makeatletter
\def\ps@pprintTitle{%
 \let\@oddhead\@empty
 \let\@evenhead\@empty
 \def\@oddfoot{}%
 \let\@evenfoot\@oddfoot}
\makeatother

\usepackage{newtxtext}

\usepackage{graphicx}
\usepackage{booktabs,colortbl}
\usepackage{xcolor}

\usepackage[leftcaption]{sidecap} 
\usepackage{booktabs}
\usepackage{lipsum}
\usepackage{subfigure}
\usepackage{graphicx}
\usepackage{bm}
\usepackage{caption}
\usepackage{mathrsfs}
\usepackage{amsmath}
\usepackage{amsfonts}
\usepackage{caption} 
\usepackage{floatrow}
\newfloatcommand{capbtabbox}{table}[][\FBwidth]
\usepackage{threeparttable,booktabs}

\usepackage{enumerate}
\usepackage{tikz}
\usepackage{listings}
\usepackage{multirow}
\usepackage{amsfonts}
\usepackage{amssymb}
\usepackage{hyperref}
\usepackage{arydshln}
\usepackage{geometry}
\newenvironment{pf}{\noindent{\textit{Proof.}  }}{\hfill $\Box$}

\def\Cov{{\rm Cov}}
 
\usepackage{float}
\hypersetup{colorlinks,%
	citecolor=blue,%
	filecolor=blue,%
	linkcolor=blue,%
	urlcolor=blue
}

\newtheorem{theorem}{Theorem}[section]

\newtheorem{rmk}[theorem]{Remark}
\newtheorem{lemma}[theorem]{Lemma}
\newtheorem{Pro}[theorem]{Proposition}

\newtheorem{problem}[theorem]{Problem}

\newtheorem{thm}[theorem]{Theorem}

\newtheorem{remark}[theorem]{Remark} 

\newcommand{\EE}{\mathbb{E}} 
\newcommand{\DD}{\mathbb{D}}

\newcommand{\RR}{\mathbb{R} }

\newcommand{\be}{\beta}

\newcommand{\al}{\alpha}

\newcommand{\om}{\omega}

\newcommand{\si}{\sigma}

\let\Section=\section
\def\section{\setcounter{equation}{0}\Section}

\def\RR{\mathbb{R} }

\def\EE{\mathbb{E}}

\def\si{{\sigma}}

\def\th{{\theta}}
\def\Th{{\Theta}}

\colorlet{tableheadcolor}{gray!10} 

\colorlet{tablerowcolor}{gray!10} 

\newcommand{\topline}{\arrayrulecolor{black}\specialrule{0.1em}{\abovetopsep}{0pt}%
            \arrayrulecolor{tableheadcolor}\specialrule{\belowrulesep}{0pt}{0pt}%
            \arrayrulecolor{black}}
 
\newcommand{\bottomlinec}{\arrayrulecolor{tablerowcolor}\specialrule{\aboverulesep}{0pt}{0pt}%
            \arrayrulecolor{black}\specialrule{\heavyrulewidth}{0pt}{\belowbottomsep}}%
\colorlet{blcolor}{gray!80}
\newcommand{\rowmidlineHR}{\arrayrulecolor{tableheadcolor}
     \specialrule{\aboverulesep}{0pt}{0pt}%
     \arrayrulecolor{black}\specialrule{\lightrulewidth}{0pt}{0pt}%
     \arrayrulecolor{tablerowcolor}\specialrule{\belowrulesep}{0pt}{0pt}%
     \arrayrulecolor{black}}

\makeatletter
\newcommand{\dotfillstretch}[1]{%
  \leavevmode
  \cleaders\hb@xt@.44em{\hss.\hss}\hskip\z@\@plus #1fill
  \kern\z@
}
\makeatother

\usepackage{natbib}
\usepackage{filecontents}

\begin{document}

\begin{frontmatter}

\title{Parameter estimation for threshold Ornstein-Uhlenbeck processes   from  discrete observations}

%
%
%

\author{Yaozhong Hu\footnote{Y.Hu is supported by an NSERC discovery grant and a startup 
fund of University of Alberta. }}  
\ead{yaozhong@ualberta.ca}
\author{Yuejuan Xi\footnote{Y. Xi  is supported by the National Natural Science Foundation of China under Grant No.  $71532001$ and $11631004$ and China Scholarship Council.}}
\ead{yjx@mail.nankai.edu.cn}

\address{Department of Math and Stat Sciences, University of Alberta at Edmonton, Canada.}
\address{School of Mathematical Sciences, Nankai University, Tianjin, China.}
\begin{abstract}
Assuming that a threshold  Ornstein-Uhlenbeck process  is observed at discrete time instants, we propose    generalized moment estimators to estimate the parameters. Our theoretical basis is    the celebrated ergodic theorem. To use this theorem we need to find the explicit form of the invariant measure.  With the sampling time step $h>0$ arbitrarily fixed, we prove the strong consistency and asymptotic normality of our  estimators as the sample size $N\to\infty$. 
\end{abstract}


\end{frontmatter}


\section{Introduction}
Let  $W=\{W(t)\}_{t\ge0}$ be a one-dimensional standard Brownian motion on a filtered probability space $ (\Omega,   \mathcal F, \mathbb P, (\mathcal F_t)_{\{t\ge 0\}} )$ and let a 
threshold Ornstein-Uhlenbeck (hereafter abbreviated as OU) process $X$ be described by   the following stochastic differential equation (SDE):
  \begin{equation}\label{e:MTOU_general}
dX_t=\sum_{i=1}^m(\beta_i-\alpha_i X_t)I(\theta_{i-1}< X_t\le \theta_i)dt+\sigma dW_t, 
\end{equation}
where  $\theta_i, i=0, 1,   \cdots, m$  with   $-\infty=\theta_0<\theta_1<\theta_2<\cdots<\theta_m=\infty$
are the so-called  thresholds;  $\beta_i\in \RR$ and $\alpha_i>0$ are the drift parameters;    $\sigma>0$ is  the diffusion parameter;      $X_0\in\RR$ is a given initial condition;    and $I(\cdot)$ denotes  the indicator function.  
 The  existence and uniqueness of the solution to the above equation \eqref{e:MTOU_general} have been known 
 \citep[e.g.][]{bass1987uniqueness}.  
 Assume that the parameters $\alpha_i$ and $\beta_i$  are unknown and assume that we can observe the state $X_t$ of the process at discrete time instants $t_k=kh$, $k=1, 2, \cdots, N$, where $h$ is an arbitrarily  given fixed time step. 
This paper aims to estimate the unknown parameters 
$\Theta=(\al_1, \cdots, \al_m, \be_1, \cdots, \be_m)$  in \eqref{e:MTOU_general} by using the obtained observations $ \left\{X_{kh}\,, k=  1, 2, \cdots, N\right\}$. 

The   models with different levels of   thresholds  
  have been widely studied and  applied  in various fields. On  one hand,
the  threshold autoregressive models  are introduced to model the nonlinearities in nonlinear time series.  
 \citet{MR717388}  found that it is more suitable  to use the 
threshold models to  describe the asymmetry in the variance-generating mechanism. 
  \citet{brockwell1991continuous}, as well as, \citet{brockwell1992continuous}   investigated   the   problems of modelling  and forecasting the continuous-time threshold  process.  
\citet{browne1995piecewise} showed that the piecewise-linear diffusion  tends to be a good approximation for some birth-and-death processes. The threshold  processes also    played an important role  in finance, we refer to  the works of  \citet{chi2017option}, \citet{decamps2006self}, \citet{jiang2018pricing}, \citet{siu2006option}, \citet{siu2016self} and references therein. 
On the other hand, 
  the threshold diffusion processes  have 
   a close tie  with the skew diffusion processes  that have  been widely treated in financial literature
   \citep[see][]{doi:10.1080/14697688.2020.1781235,gairat2017density,wang2015skew,zhuo2017efficient,zhuo2017simple}.

While the threshold models are applied, an important problem is to estimate the parameters $\Theta$ through the  available historical 
data. There have been  some approaches to estimate the parameters for threshold diffusion processes such as least squares estimation, likelihood estimation, and Bayesian estimation. We refer the readers  to \citet{brockwell2007continuous},  \citet{chan1993consistency},   \citet{kutoyants2012identification},  \citet{doi:10.1111/sjos.12417},  and \citet{stramer2007bayesian}. Let us also mention that in \citet{su2015quasi, su2017testing}, the authors  proposed  the  novel quasi-likelihood estimators and test.
Within the above mentioned  estimation methods, the observations are supposed to be obtained  continuously. 
Since real data are usually collected at  discrete  time instants,  it is 
necessary to estimate the parameters when only discrete observations are available.  
To our best knowledge, the problem to  estimate parameters  for  a continuous-time threshold diffusion  processes based on discrete observations is   under-explored.  
 
One situation in the discrete-time observations is that one has  the  high-frequency data, which means that in our observations
$\left\{X_{kh}, k=1, 2, \cdots, N\right\}$, we have $h$ depends on $N$,  $h\rightarrow 0$,  and $Nh\rightarrow \infty$. In this case it is possible to approximate the (stochastic) integral by its ``Riemann-It\^o" sum  to modify the continuous-time estimators to the discrete ones.  

In reality,  the continuous   or  high-frequency observations are   usually  impossible or very costly that we cannot
have the luxury to collect such large amount of data. As a consequence, 
the time step $h$ must be allowed to be an arbitrarily  fixed constant. Hence, 
we cannot borrow   methods that are  only valid 
for continuous-time observations or  for high-frequency data.  The present work proposes a completely different approach to address this problem. Our approach  is motivated by the previous works     of 
the construction of  the  estimators:   the ergodic type estimators for the OU process driven by fractional Brownian motion \citep[e.g.][]{hu2013parameter}; the ergodic type estimators for  the reflected  OU process driven by standard  Brownian motions
 \citep[e.g.][]{hu2015parameter}; and  the ergodic type estimators for the OU process driven     by stable L\'evy motions \citep[e.g.][]{cheng2020generalized}.

Similar to the above mentioned papers, we use the ergodic theorem to obtain the generalized 
moment estimators for the parameters.  To this end, we need first to prove the ergodic theorem for our threshold diffusion 
\eqref{e:MTOU_general}. 
Namely, we need to prove  that there is a probability 
density function $\psi(x)$  such that
\begin{equation*}
\lim_{N\rightarrow \infty}\frac{1}{N} \sum_{k=1}^N f(X_{t_k})=\int_{\RR} f(x) \psi(x) dx 
\end{equation*}
and we also need to find the explicit form of the probability density $\psi(x)$. 
This is done in Section \ref{sec:pre}.  After obtaining the    explicit dependence of the probability density  on the parameters we let 
\begin{equation}
 \frac{1}{N}\sum_{k=1}^N f_i(X_{t_k})=\int_{\RR} f_i(x) \psi(x) dx \label{e.1.2}
\end{equation}
for    different appropriately chosen functions $f_i$ to obtain a  suitable system of  algebraic equations for the parameters. 
In Equation \eqref{e:MTOU_general} there are $3m$ unknown parameters
$ \al_1, \cdots, \al_m, \be_1, \cdots, \be_m, \th_1, \cdots, \th_{m-1}, \si$. Presumably, we can choose $3m$ different functions $f$ so that we obtain  a system of 
$3m$ equations for the $3m$ unknowns. However,  some parameters are coupled with each other and cannot be separated. For example, from Remark \ref{r.3.1},
when $m=2$, $ \th_0=-\infty $, $\th_1=0$, $\theta_2=\infty$, 
$\be_1=\be_2=0$, we see that if 
$(\frac{\al_1}{\si^2}, \frac{\al_2}{\si^2})$ remains the same, then the invariant probability density $\psi_1$ remains the same function. So, even in this simplest case we cannot expect to use \eqref{e.1.2} to estimate $\al_1$, $\al_2$, and $\si$ simultaneously.   To avoid this  identifiability problem   in this paper we focus on the estimation of the parameters $\Theta$ assuming $\th_1, \cdots, \th_{m-1}, \si$ are known.   Furthermore, to better convey   our idea,  we focus on the case that $m=2$, $\th_0=-\infty$, $\th_2=\infty$, and the parameters $\th=\th_1$ and $\si$
are known.   This means that we shall focus on the following equation:
\begin{equation}\label{e:MTOU}
dX_t= (\beta_1-\alpha_1 X_t)I(  X_t\le \theta )dt+
(\beta_2-\alpha_2 X_t)I(  X_t> \theta )dt+
\sigma dW_t\,,
\end{equation}
where  $\theta\in\RR$, $\beta_1$, $\beta_2\in\RR$, $\alpha_1$, $\alpha_2\in(0,\infty)$, and $\sigma\in(0,\infty)$.  However, it should be mentioned that if $\si$ and $\theta$ are unknown, we may  assume that the data are collected from the high-frequency type. In this case, $\si$ and $\theta$ can be estimated in the manners of \citet{kutoyants2012identification} and  \citet{su2015quasi}, respectively. Now that we have four parameters $\Th=(\al_1, \al_2, \be_1, \be_2)$, so we only need to choose four different $f$ to obtain
a system of  four   equations. However, since the invariant probability density
 $\psi$  depends on the parameters in a very complex  way 
it is hard to know whether  the solution exists  (locally and globally) uniquely. One of the major contributions of this work 
is to  appropriately use the conditional moments   so that we can obtain  some manageable  equations. This will be carried 
out in Section \ref{sec:S}. We briefly   summarize our efforts
in that section  as follows. 
\begin{enumerate} 
\item[(1)] In Section \ref{sub:I}, we assume $\be_1=\be_2 =\th=0$. The conditional moments are introduced to obtain the explicit 
generalized moment estimators for $\al_1$ and $\al_2$. Furthermore, the strong consistency and asymptotic normality of the estimators are obtained. 
\item[(2)] In Section \ref{sub:II}, we assume that $\be_1=\be_2 =0$ whereas $\th$ is known but is not equal to $0$.   In this case, we can obtain  two uncoupled algebraic equations  for the two parameters $\al_1$ and $\al_2$ by conditional moments. 
Each of these equations will be shown to have a globally unique solution, yielding the generalized moment estimators for $\al_1$ and $\al_2$, although not explicitly.  The strong consistency and asymptotic normality of the estimators are obtained. 
\item[(3)] In Section \ref{sub:III},  we further assume that $\th$ is known but is equal to not $0$ and we want to estimate all the four parameters $(\al_1, \al_2, \be_1, \be_2)$. We use  the conditional moments to invert the four equations into two uncoupled systems of  equations to obtain the generalized estimators for $\alpha_1$, $\alpha_2$, $\beta_1$, and $\beta_2$.
The Jacobians (which are independent of data) of the two systems are computed, whose 
non-degeneracy implies that both systems have unique local solutions. To seek an answer for  global uniqueness we reduce the problem to a simpler  one  of  finding the zeros of two functions, both  of a single variable. If the derivatives (now involving observation data) of such functions are nonzero, then   the global uniqueness holds  by the mean value theorem.  
\end{enumerate} 
  In our cases (2) and (3)  the explicit solution to the system of algebraic equations is still hard to obtain. But there are many standard  methods, such as the Newton-Raphson  iteration method. It is  available to solve the nonlinear system in Matlab and Mathematica by the built-in functions ``fsolve'' and ``FindRoot'', respectively. 
%
In Section \ref{sec:N},    some numerical experiments are provided to show the efficiency of  our estimation approach. Section \ref{sec:C} concludes this paper.
 
\section{Ergodicity and invariant density}
\label{sec:pre}
Before proceeding to construct our  estimators, we need some stationary and ergodic properties of the threshold diffusion process described by \eqref{e:MTOU_general}. The following  proposition  is adopted from \citet{brockwell1991continuous}, \citet{brockwell1992continuous}, 
and \citet{browne1995piecewise}. 
\begin{Pro}\label{pro:density}
Suppose that $\sigma>0$.  Then the process defined by \eqref{e:MTOU_general} has a stationary distribution if and only if 
\begin{equation*}
\lim_{x\to-\infty}(-\alpha_1 x^2+2\beta_1x)<0\,,\quad \lim_{x\to\infty}(-\alpha_m x^2+2\beta_mx)<0\;.
\end{equation*}
Furthermore, the stationary density is given by
\begin{equation*}
\psi(x)=\sum_{i=1}^mk_i\exp\left(\frac{-\alpha_ix^2+2\beta_ix}{\sigma^2}\right)I(\theta_{i-1}<x\le \theta_i),
\label{e.2.1} 
\end{equation*}
where $k_i$ are uniquely determined by 
the system of $m$ equations: 
\begin{equation*}
\int_{-\infty}^\infty\psi(x)dx=1\,,\quad {\rm and}\quad \psi(\theta_i-)=\psi(\theta_i+)\,,\quad i=1,2,\ldots, m-1 \,.
\end{equation*} 
\end{Pro}
\begin{rmk}\label{rmk:0}   The constants $k_i, i=1, 2, \cdots, m$ depends on the parameters in
the equation \eqref{e:MTOU_general}.  This is one of the main reasons to make the analysis of the system of algebraic equations  sophisticated. 
\end{rmk}
Although the stationary density function $\psi(\cdot)$ is not Gaussian, it is a mixture of   Gaussian densities and has finite moments of all orders. Moreover, if the threshold OU process $X$ is stationary, it is also geometrically ergodic  \citep[see][]{stramer1996existence}. The following lemma describes the stochastic stability of threshold OU processes and plays a crucial role in our estimation approach.
\begin{lemma}\label{le:ergodic}
The $h$-skeleton sampled chain $\{X_{kh}:k\ge0\}$ which comes from  the process $X$ defined by  \eqref{e:MTOU_general} is ergodic, namely, the following   ergodic identity holds: for any $X_0\in\mathcal S:=\RR$ and  for  any $f\in L_1(\RR, \psi(x) dx)$,
\begin{equation*}
\lim_{N\to\infty}\frac{1}{N}\sum_{k=1}^N f(X_{kh})=\mathbb E[f(X_\infty)]=\int_{\mathbb R} f(x)\psi(x)dx, ~a.s.
\end{equation*}
\end{lemma}
\begin{pf} 
It suffices to show that the process X is bounded in probability on average and is a $T$-process \citep[see][Theorem 8.1]{meyn1993stability}. We note that  the threshold diffusion process $X$ is a $\varphi$-irreducible $T$-process, where $\varphi$ is a Lebesgue measure  \citep[see][]{stramer1996existence}. Moreover, since for  $i=1,2$,
\begin{equation*}
\lim_{|x|\to\infty}\left(-\alpha_i x^2+2\beta_i x \right)<0,
\end{equation*}
 we have from \citet[Theorem 5.1]{stramer1996existence} that  $X$  is a positive Harris recurrent process.
Finally, by virtue of  \citet[Theroem 3.2(ii)] {meyn1993stability}, we   conclude that $X$  is bounded in probability on average.
\end{pf}


 Using the same definitions as that in \citet{karlin1981second},  the scale density function $s(x)$, scale measure $S(x)$, and speed density function  $m(x)$ are given by
\begin{equation*}
s(x)=\left\{
     \begin{aligned}
        &c_1\exp\left( -\frac{2\beta_1x}{\sigma^2}+\frac{\alpha_1x^2}{\sigma^2}\right) \;, & x\le\theta,\\
        &c_2\exp\left( -\frac{2\beta_2x}{\sigma^2}+\frac{\alpha_2x^2}{\sigma^2}\right)\;, &x>\theta,
     \end{aligned}
     \right.
\end{equation*}
\begin{equation*}
 S(x)=\int_{-\infty} ^x s(y)dy, \quad m(x)=\frac{2}{s(x)\sigma^2}, 
\end{equation*}
where $c_1=\exp\left(-\frac{2\beta_2\theta}{\sigma^2}+\frac{\alpha_2\theta^2}{\sigma^2}\right)$ and $c_2=\exp\left(-\frac{2\beta_1\theta}{\sigma^2}+\frac{\alpha_1\theta^2}{\sigma^2}\right)$. For $i=1,2$, let 
\begin{equation*}
\widetilde z_i=\frac{\sqrt{2\alpha_i}}{\sigma}\left(\theta-\frac{\beta_i}{\alpha_i}\right),  \quad b_i=\frac{\beta_i^2}{\sigma^2\alpha_i}\;.
\end{equation*}
Then the coefficients $k_1$ and $k_2$ of $\psi(x)$ are given by
\begin{align}\label{e:k1}
k_1&=\frac{1}{\sigma\sqrt \pi}\frac{\phi(\widetilde z_2)}{\phi(\widetilde z_2)e^{b_1}\Phi(\widetilde z_1)/\sqrt{\alpha_1}+{\phi(\widetilde z_1)e^{b_1}\Phi(-\widetilde z_2)/\sqrt{\alpha_2}}
}\;,\\\label{e:k2}
k_2&=\frac{1}{\sigma\sqrt \pi}\frac{\phi(\widetilde z_1)}{\phi(\widetilde z_2)e^{b_2}\Phi(\widetilde z_1)/\sqrt{\alpha_1}+{\phi(\widetilde z_1)e^{b_2}\Phi(-\widetilde z_2)/\sqrt{\alpha_2}}
}\;,
\end{align}
where $\phi(x):=\frac{1}{\sqrt{2\pi}}e^{-x^2/2}$ is the normal density, $\Phi(x)=\frac{1}{\sqrt{2\pi}}\int_{-\infty}^xe^{-\frac{y^2}{2}}dy$ is the the standard normal distribution function.
Although the SDE \eqref{e:MTOU} has no explicit solution, we can derive the  spectral expansion of its transition density,   see the proof in \ref{A:spectral}.

\begin{Pro}\label{Pro:spec}
For $i=1,2$, set
\begin{align*}
z_i&=\frac{\sqrt{2\alpha_i}}{\sigma}\left(x-\frac{\beta_i}{\alpha_i}\right),  \quad \nu_i=\frac{\lambda}{\alpha_i} \;,\\
\widetilde z_i&=\frac{\sqrt{2\alpha_i}}{\sigma}\left(\theta-\frac{\beta_i}{\alpha_i}\right), \quad \varrho=\frac{2\beta_1\theta+2\beta_2\theta-\alpha_1\theta^2-\alpha_2\theta^2 }{\sigma^2}.
\end{align*}
Let   $D_v(z)$ and $H_v(z)$ denote  the parabolic cylinder function and Hermite function respectively \citep[see][]{MR0240343,MR0174795}. 
Let  $0\le\lambda_1<\lambda_2<\cdots<\lambda_n\to\infty$ as $n\to\infty$ be  the simple discrete zeros of the Wronskian equation: 
\begin{equation}\label{e:w}
\omega(\lambda)=\exp(\varrho)2^{1-\frac{\nu_1+\nu_2}{2}}\sigma^{-1}\left[\nu_2\sqrt{\alpha_2} H_{\nu_1}(-\frac{\widetilde z_1}{\sqrt2})H_{\nu_2-1} (\frac{\widetilde z_2}{\sqrt 2}) +\nu_1 \sqrt{\alpha_1} H_{\nu_2}(\frac{\widetilde z_2}{\sqrt2})H_{\nu_1-1} (-\frac{\widetilde z_1}{\sqrt 2})\right] =0\,. 
\end{equation}
Denote 
\begin{equation}\label{e:varphi}
\varphi_n(x)=\left\{
     \begin{aligned}
        &\sqrt{\frac{\eta(\theta,\lambda_n)}{\omega^\prime(\lambda_n)\xi(\theta,\lambda_n)}}\xi(x,\lambda_n)  \;, & x\le\theta,\\
        &\sign(\xi(\theta,\lambda_n)\eta(\theta,\lambda_n))\sqrt{\frac{\xi(\theta,\lambda_n)}{\omega^\prime(\lambda)\eta(\theta,\lambda_n)}}\eta(x,\lambda_n)  \;, &x>\theta,
     \end{aligned}
     \right.
\end{equation}
with 
\begin{equation*}
\xi(x,\lambda)=\exp\left( z_1^2/4\right)D_{\nu_1}(-z_1), \quad \eta(x,\lambda)=\exp\left( z_2^2/4\right)D_{\nu_2}(z_2)\,. 
\end{equation*}
  Then, 
the spectral expansion of the transition  density of $X$  (defined from
$ \mathbb P(X_t\in A|X_0=x)=\int_A p_t(x,y)dy$ for any Borel set
$A$  of
$\RR$)    is given by
\begin{equation*}\label{density}
p_t(x,y)=m(y)\sum_{n=1}^{\infty}\exp(-\lambda_n t)\varphi_n(x)\varphi_n(y)\,. 
\end{equation*}
\end{Pro}

\section{Estimate  $\alpha_i$ and $\beta_i$}
\label{sec:S}
In this section  we attempt to construct   generalized  moment estimators  for the parameters $\alpha=(\alpha_1,\alpha_2)^T$ and $\beta=(\beta_1,\beta_2)^T$, where $T$ denotes the transpose of a vector, and  to study their   strong consistency and asymptotic normality.  We classify our study  into several cases according to the drift parameters.

\subsection{Case I: Estimate $\alpha_i$ for known $\beta_i=0$ and $\theta=0$ }\label{sub:I}
Here we consider the case $\beta_i=0$, $i=1,2$  and $\th=0$.   In this case the equation 
becomes
\begin{equation}\label{e:0OU}
dX_t =-\alpha_1 X_t I( X_t\le 0)dt  -\alpha_2 X_t I( X_t> 0)dt+\sigma dW_t\,. 
\end{equation}
Then the stationary density of $X$ is given by  
\begin{equation}
\psi_1(x)=\frac{2 \sqrt{\alpha_1\alpha_2}}{\sqrt{\pi}(\sqrt{\alpha_1}+\sqrt{\alpha_2}) \sigma }\left[
\exp\left(-\frac{\alpha_1x^2}{\sigma^2}\right)  I(x\le0)+ \exp\left(-\frac{\alpha_2x^2}{\sigma^2}\right)  I(x>0)
\right].\label{e.3.6} 
\end{equation}
\begin{remark}\label{r.3.1}  It is easily observed that  
$\psi_1(x)$ depends only on  $\frac{\al_1}{\si^2}$ 
and $\frac{\al_2}{\si^2}$. 
\end{remark}
From this identity we have 
\begin{Pro}\label{Prop:Xestimate}
Let $X_\infty=\lim_{t\to\infty}X_t$ and define
\begin{equation*}\label{e:LR}
L_n=\mathbb E\left[(-X_\infty)^{n} I(X_\infty\le0)\right], \quad R_n=\mathbb E\left[X_\infty^{n} I(X_\infty>0)\right].
\end{equation*}
Then for any real number $n>0$,
\begin{align}
L_n&= \frac{\sigma^{n}\sqrt{\alpha_1\alpha_2}}{\alpha_1^{(n+1)/2}\sqrt{\pi}(\sqrt{\alpha_1}+\sqrt{\alpha_2})}\Gamma\left(\frac{n+1}{2}\right)\;, \label{e.3.8} \\
R_n&=\frac{\sigma^{n}\sqrt{\alpha_1\alpha_2}}{\alpha_2^{(n+1)/2}\sqrt{\pi}(\sqrt{\alpha_1}+\sqrt{\alpha_2})}\Gamma\left(\frac{n+1}{2}\right)\;,\label{e.3.9} 
\end{align}
where $\Gamma(\cdot)$ denotes the Gamma function $\Gamma(\alpha)=\int_0^\infty x^{\alpha-1}e^{-x}dx $.
\end{Pro}

From the above expressions \eqref{e.3.8}-\eqref{e.3.9} and by some elementary calculations, we  can represent the parameters $\alpha_1$ and $ 
\al_2$ in terms of $L_n$ and $R_n$  as
\begin{align}\label{e:Ealpha}
\alpha_1&=\left\{\frac{ \sigma^{n}\Gamma\left(\frac{n+1}{2}\right)}{\sqrt{\pi}L_n\left[ \left( \frac{R_n}{L_n}\right)^{\frac{1}{n+1}}+1\right]} \right\}^{\frac{2}{n}},\\\label{e:Ebeta}
\alpha_2&=\left\{\frac{\sigma^{n}\Gamma\left(\frac{n+1}{2}\right)}{\sqrt{\pi}R_n\left[ \left(  \frac{L_n}{R_n}\right)^{\frac{1}{n+1}}+1\right]} \right\}^{\frac{2}{n}}.
\end{align}
Since $ L_n>0$  and $R_n>0$, $\al_1$ and $\al_2$   
 are well-defined by \eqref{e:Ealpha} and \eqref{e:Ebeta}.

Setting
\begin{equation*}
\widehat L_{n,N}=\frac{1}{N}\sum_{k=1}^N(-X_{kh})^{n}I(X_{kh}\le0),\quad
\widehat R_{n,N}=\frac{1}{N}\sum_{k=1}^N(X_{kh})^{n}I(X_{kh}>0),
\end{equation*}
 we naturally construct  the generalized  moment estimators for $\al_1, \al_2$ as  follows: 
\begin{align}\label{e:alpha}
\widehat\alpha_{1,n,N}&=\left\{\frac{ \sigma^{n}\Gamma\left(\frac{n+1}{2}\right)}{\sqrt{\pi}\widehat L_{n,N}\left[ \left( \frac{\widehat R_{n,N}}{\widehat L_{n,N}}\right)^{\frac{1}{n+1}}+1\right]} \right\}^{\frac{2}{n}},\\\label{e:beta}
\widehat\alpha_{2,n,N}&=\left\{\frac{\sigma^{n}\Gamma\left(\frac{n+1}{2}\right)}{\sqrt{\pi}\widehat R_{n,N}\left[ \left( \frac{\widehat L_{n,N}}{\widehat R_{n,N}}\right)^{\frac{1}{n+1}}+1\right]} \right\}^{\frac{2}{n}}.
\end{align}
We will show the strong consistency and asymptotic normality of the estimators $\widehat\alpha_{1,n,N}$ and $\widehat\alpha_{2,n,N}$ of $\alpha_{1}$ and $\alpha_{2}$ in the following theorems.   
\begin{rmk}
Although the expectation of $(-X_\infty)^n I(X_\infty\le \theta)$ (or $X_\infty^n I(X_\infty>\theta)$)
 is not the $n$-th  order moment in the conventional sense, it captures   sufficient information about  the parameters and  the motivation of the estimation scheme in this paper stems from the generalized 
  moment estimation. For this reason, we still use the term of ``generalized  moment estimators".
\end{rmk}
\begin{thm}\label{thm:1}
Fix any real number $n>0$ and fix any time step size 
$h>0$. Then $\widehat\alpha_{1,n,N}\to\alpha_1$ and $\widehat\alpha_{2,n,N}\to\alpha_2$ almost surely as $N\to\infty$, where  $\widehat\alpha_{1,n,N}$, $\widehat\alpha_{2,n,N}$ 
are  defined by 
\eqref{e:alpha} and \eqref{e:beta} 
 respectively. 
\end{thm}
\begin{pf} The   straightforward applications of Lemma \ref{le:ergodic} to  $f_1(x)=(-x)^nI(x\le 0)$ and $f_2(x)=x^nI(x>0)$  yield 
\begin{equation*}
\lim_{N\to\infty} \widehat L_{n,N}= L_n> 0,\quad \lim_{N\to\infty}\widehat R_{n,N}=R_n>0,~a.s. 
\end{equation*}
which  imply the theorem by \eqref{e:Ealpha}-\eqref{e:beta}. 
\end{pf}

Next, we study the central limit theorem (CLT) for the estimators. In comparison to Theorem 2 in \citet{hu2015parameter},  we shall  discuss the joint asymptotic normality of the estimators. Before stating our theorem we need the following notations. Denote 
\begin{equation*}
g_{1n}(x)=(-x)^nI(x\le0), \quad g_{2n} (x)=x^nI(x>0) \,. 
\end{equation*}
Let $\widetilde X_0$  be a random variable with probability density function $\psi_1$
given by \eqref{e.3.6}, independent of the Brownian motion and let
$\widetilde X_t$ be the solution to \eqref{e:0OU} with initial condition $\widetilde X_0$. 
From \citet[Theorem 17.0.1]{meyn2012markov}, we get that
\begin{equation}\label{e:sigma}
\sigma_{ij}^n:=\Cov ( g_{in} (\widetilde X_0),g_{jn}(\widetilde X_0))+\sum_{k=1}^\infty\left[ \Cov ( g_{in} (\widetilde X_0),g_{jn} (
\widetilde X_{kh}))+\Cov (g_{jn} (\widetilde X_0), g_{in}  (\widetilde X_{kh})) \right]\,, 
\end{equation}
where $i,j=1,2$, are well defined and  are  given by \eqref{e.B.2} with $\th=0$.    
Let   $G_{i,n}$, $i=1,2$ be   defined on $\RR^2$ by 
\begin{equation*}
G_{1,n}(x,y)=\left\{\frac{ \sigma^{n}\Gamma\left(\frac{n+1}{2}\right)}{\sqrt{\pi}x\left[ \left( \frac{y}{x}\right)^{\frac{1}{n+1}}+1\right]} \right\}^{\frac{2}{n}},\quad
G_{2,n}(x,y)=\left\{\frac{\sigma^{n}\Gamma\left(\frac{n+1}{2}\right)}{\sqrt{\pi}y\left[ \left( \frac{x}{y}\right)^{\frac{1}{n+1}}+1\right]} \right\}^{\frac{2}{n}}
\end{equation*}
which are  the functions corresponding to \eqref{e:alpha} and \eqref{e:beta}.  Denote $G_n=(G_{1,n},G_{2,n}): \RR^2\to \RR ^2$.

\begin{figure} [t]
\centering 
\subfigure[$\widehat\alpha_{1,n,N}$ against $n$ ($\al_1=0.02$)]{\label{fig:a1}
\includegraphics[width=0.4\textwidth]{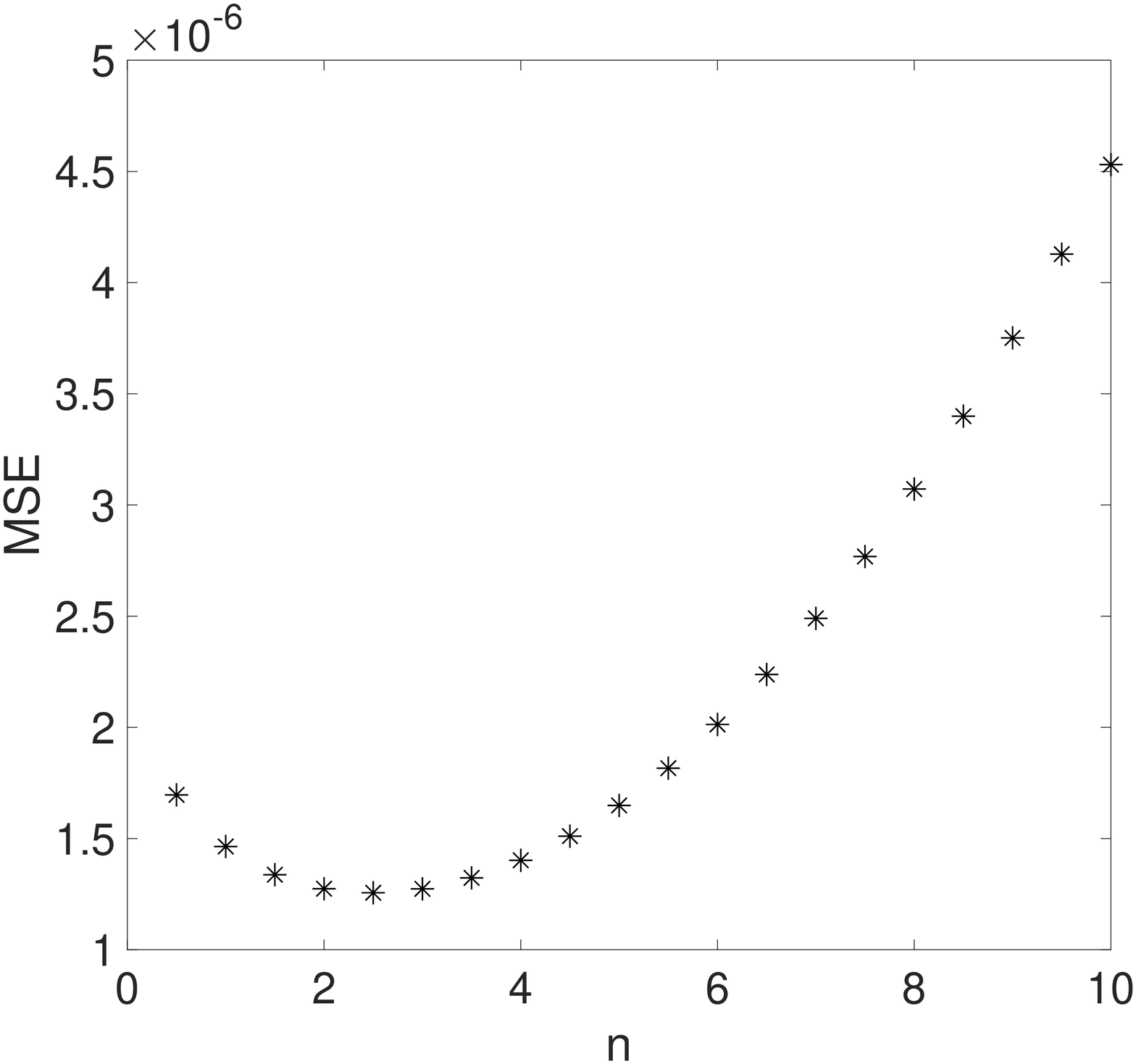}}
\subfigure[ $\widehat\alpha_{2,n,N}$ against $n$ ($\al_1=0.05$)]{\label{fig:b1}
\includegraphics[width=0.4\textwidth]{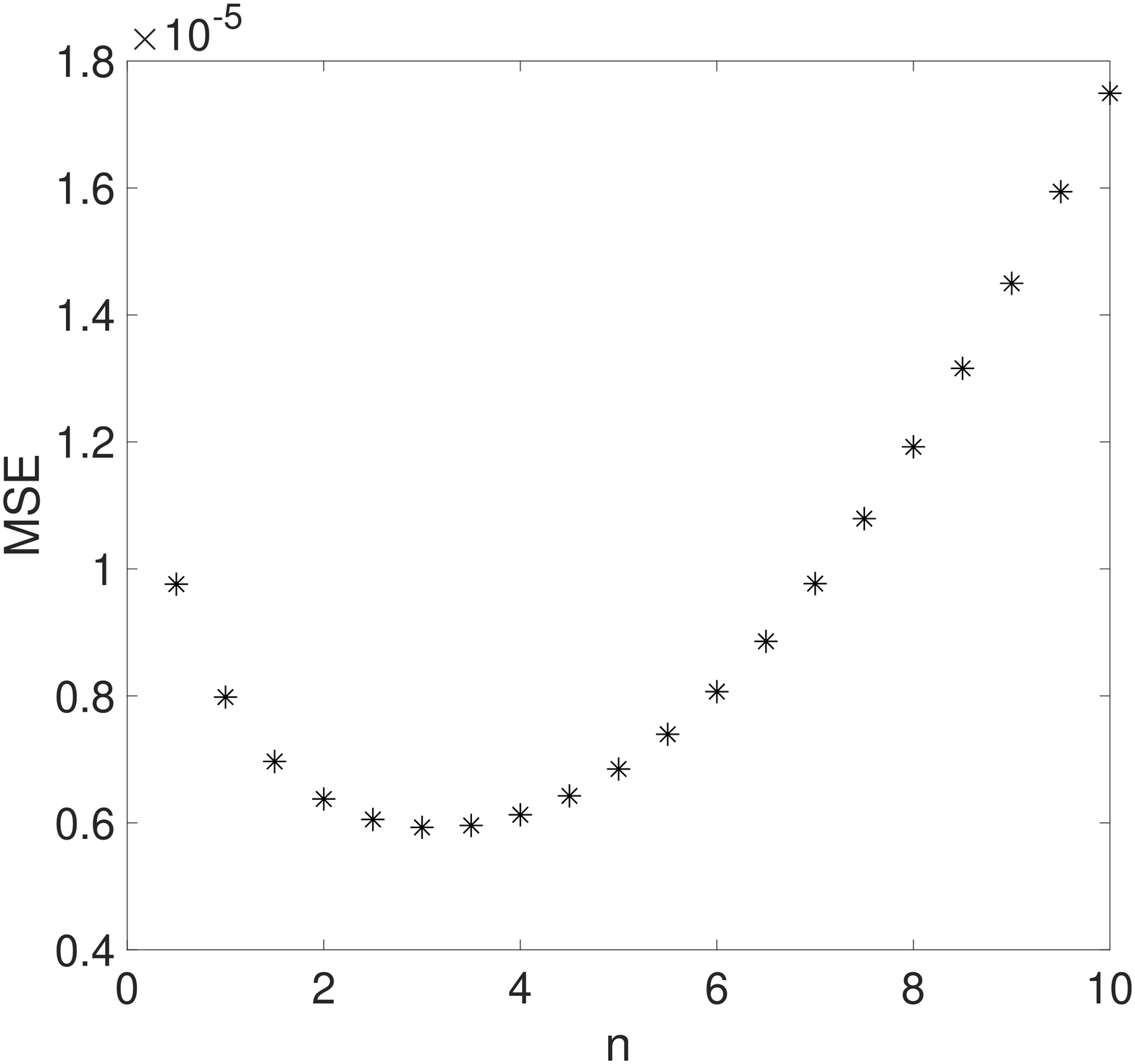}}
\subfigure[$\widehat\alpha_{1,n,N}$ against $n$ ($\al_1=0.1$)]{\label{fig:b2}
\includegraphics[width=0.4\textwidth]{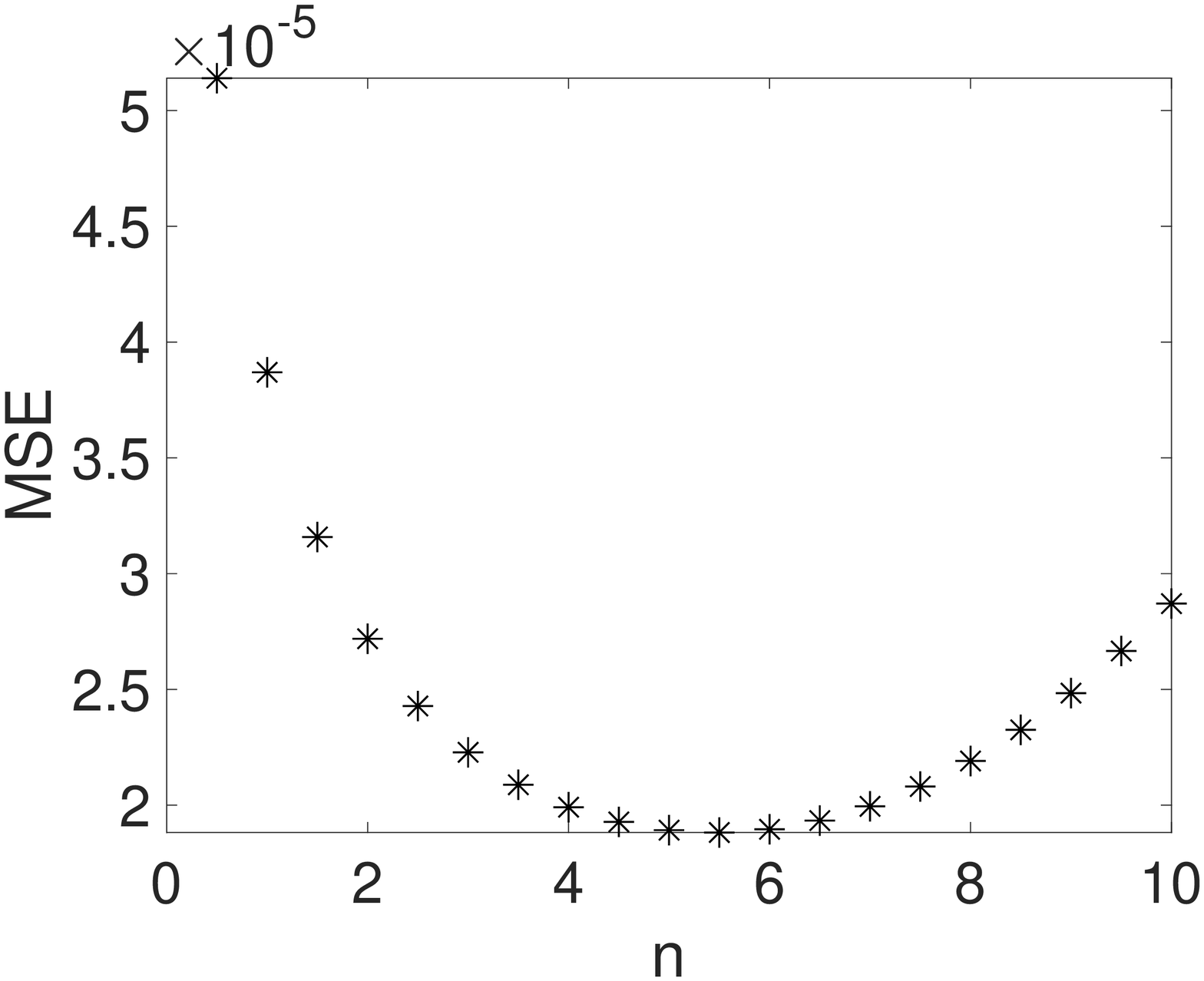}}
\subfigure[ $\widehat\alpha_{2,n,N}$ against $n$ ($\al_1=0.5$)]{\label{fig:b2}
\includegraphics[width=0.4\textwidth]{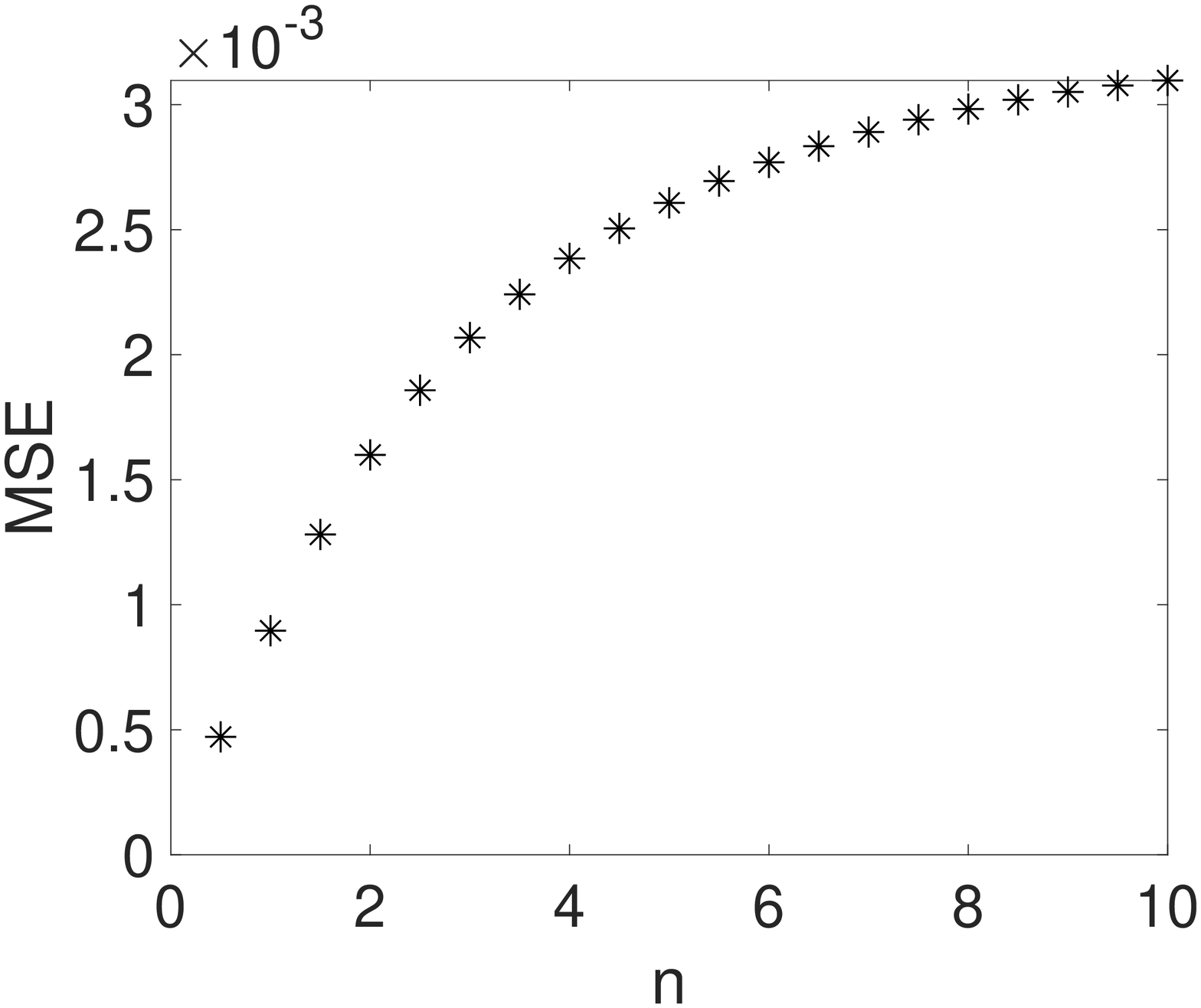}}
\caption{MSE of  $\widehat\alpha_{1,n,N}$ and  $\widehat\alpha_{2,n,N}$. }
\label{fig:MSE}
\end{figure}

Now we can state our main result of this subsection. 
\begin{thm}\label{asy:1}
Fix an arbitrary $h>0$.  Denote $\alpha=(\alpha_1,\alpha_2)^T$ and $\widehat\alpha_{n, N}=(\widehat\alpha_{1,n,N},\widehat\alpha_{2,n,N})^T$,  where  $\widehat\alpha_{1,n,N}$, $\widehat\alpha_{2,n,N}$ 
are  defined by
\eqref{e:alpha} and \eqref{e:beta} 
 respectively. Then as $N\to\infty$,
\begin{align*}
\sqrt N\left(\widehat\alpha_{n, N} -\alpha \right)&\Rightarrow \mathbf N\left(0, \nabla G_n(L_n, R_n) \cdot\Sigma_n\cdot \nabla G_n(L_n,R_n)^T\right)\,, 
\end{align*}
where   the symbol ``$\Rightarrow$'' denotes convergence in distribution,  $\mathbf N(\mu,\Sigma)$ stands for the normal  random vector  with mean $\mu$ and  variance $\Sigma$, and $ \Sigma_n:=(\sigma_{ij}^n)_{1\le i,j\le 2}$  with  $\sigma_{ij}^n$ being  defined by \eqref{e:sigma} or equivalently by \eqref{e.B.2} with $\th=0$\,.  
\end{thm}
\begin{pf}
The proof is carried out  in two steps. First, we establish the bivariate CLT for $(\widehat L_{n,N}, \widehat R_{n,N})^T$, then we employ the bivariate delta method. Recall that $\{X_{kh}\}$ is a positive Harris chain with invariant probability $\psi$ (Lemma \ref{le:ergodic}) and  is $V$-uniformly ergodic with a function $V(x)=x^{2m}+1$ or $V(x)=e^{x^{2m}}+1$ 
 \citep[see][Theorem 5.1]{stramer1996existence}. That is to say, there exist $R\in(0,\infty)$ and $\rho\in(0,1)$ such that for all $x\in \mathbb R$,
\begin{equation*}
||P^n(x,\cdot)-\psi_1 ||_V\le RV(x)\rho^n,
\end{equation*}
where $V$-norm  $||\nu||_V:=\sup_{g:g\le V}|\nu(g)|$, $\nu$ is any signed measure \citep[see][Page 334]{meyn2012markov}, and $P^n(x, B):=P_{nh}(x, B):=\mathbb P(X_{nh}\in B)$ is an $n$-step transition probability function of the sampled chain $\{X_{kh}\}_{k\ge0}$ from the initial point $x$ to set $B$. Then from \citet[Theorem 17.0.1]{meyn2012markov}, for any $(a_1, a_2)\in \mathbb R^2$,  letting $A(x)=a_1x^nI(x\le0)+a_2x^nI(x>0)$, we know that
\begin{equation*}
\sqrt N (a_1 \widehat L_{n,N} +a_2 \widehat R_{n,N})=\frac{1}{\sqrt N} \sum_{k=1}^N A(X_{kh})=:\frac{1}{\sqrt N}S_n(A)
\end{equation*}
converges  to some normal random variable $Z_{a_1, a_2}$  in the sense of distribution. 
By the Cram\'er-Wold device 
 we know that
$\sqrt N(\widehat L_{n,N} , \widehat  R_{n,N})^T$ converges jointly to a (two-dimensional) normal vector.  Moreover, in view of the multivariable Markov chain CLT \citep[Section 1.8.1]{brooks2011handbook}, we have 
\begin{equation*}
\sqrt N((\widehat L_{n,N} , \widehat  R_{n,N})^T-( L_{n,N} ,   R_{n,N})^T )\Rightarrow \mathbf N(0, \Sigma_n)\,, 
\end{equation*}
where $\sigma_{ij}^n$ is  defined by \eqref{e:sigma}.   
Let us recall    the sufficient conditions of the multivariate delta method  \citep[see][]{MR1652247}: all partial derivatives $\partial G_j(x,y)/\partial x$ and  $\partial G_j(x,y)/\partial y$  exist for $(x,y)$ in a neighborhood of $(L_n, R_n)$ (notice that $L_n>0$ and $R_n>0$) and are continuous at $(L_n,R_n)$. It is clear that the conditions are justified.
From the multivariate delta method, the following desired result follows 
\begin{equation*}
\sqrt N \left( G_n(\widehat L_{n,N}, \widehat R_{n,N})-G_n(L_n,R_n)\right)\Rightarrow \mathbf{N}\left(0, \nabla G_n(L_n, R_n) \cdot\Sigma_n\cdot \nabla G_n(L_n,R_n)^T\right).
\end{equation*}
Hence, we complete the proof.
\end{pf} 
\begin{rmk}
The asymptotic variance is given by $ \nabla G_n(L_n, R_n) \cdot\Sigma_n\cdot \nabla G_n(L_n,R_n)^T$. Our numerical experiments show that the estimators perform better when $\alpha_1$ and $\alpha_2$ are smaller in terms of mean squared error (MSE), see Figure \ref{fig:MSE}, where we set $\sigma=1$, $h=0.5$, $N=100,000$.  From  Figure \ref{fig:MSE} , we also see   that the best estimators is to choose the moment $n$ to be between $2$ and $4$.
\end{rmk}
%
 \subsection{Case II: Estimate $\alpha_i$ for known $\beta_i=0$ and $\theta\neq0$} \label{sub:II}
Now we consider the case    $\theta\neq0$,  $\beta_i=0$,   $i=1,2$. Recall the explicit expression for  the stationary density we
obtained in   Section \ref{sec:pre}: 
\begin{equation*}
\psi_2(x)=k_1\exp\left(-\frac{\alpha_1x^2}{\sigma^2}\right)I(x\le \theta)+k_2\exp\left(-\frac{\alpha_2x^2}{\sigma^2}\right)I(x> \theta),
\end{equation*}
where $k_1$ and $k_2$ are determined by $\psi_2(\theta-)=\psi_2(\theta+)$ and $\int_{-\infty}^\infty\psi_2(x)dx=1$.
The constants $k_1$ and $k_2$ are  complicated functions of 
the unknown parameters $\al_1$ and $\al_2$. We shall use the technique of conditional moments to get rid of them. Since the stationary distribution of $X$ is Gaussian conditioned to stay in the interval  $(-\infty,\theta)$ or the interval 
  $(\theta,\infty)$, we shall focus on the conditional moments of $X_\infty$. Some elementary calculations give
\begin{equation*}
    \left\{
     \begin{aligned}
        &\EE [X_\infty|X_\infty\le\theta]=\frac{\EE [X_\infty I(X_\infty\le\theta)]}{\EE [I(X_\infty\le\theta)]}=-\frac{\sigma}{\sqrt{2\alpha_1}}\frac{\phi( -\sqrt{2\alpha_1}\theta/\sigma)}{1-\Phi(-\sqrt{2\alpha_1}\theta/\sigma)},\\
        &\EE [X_\infty|X_\infty>\theta]=\frac{\EE [X_\infty I(X_\infty>\theta)]}{\EE [I(X_\infty>\theta)]}=\frac{\sigma}{\sqrt{2\alpha_2}}\frac{\phi( \sqrt{2\alpha_2}\theta/\sigma)}{1-\Phi(\sqrt{2\alpha_2}\theta/\sigma)}\,. 
     \end{aligned}
     \right.\label{e.3.12}
\end{equation*}
For simplicity of notations,
we set 
\begin{equation}\label{LR}
\widehat L_{n,N}^\theta=\frac{1}{N}\sum_{k=1}^NX_{kh}^{n}I(X_{kh}\le\theta),\quad
\widehat R_{n,N}^\theta=\frac{1}{N}\sum_{k=1}^NX_{kh}^{n}I(X_{kh}>\theta),
\end{equation}
\begin{equation}
 L_{n}^\theta=\EE [X_\infty^n I(X_\infty\le\theta)],\quad
 R_{n}^\theta=\EE [X_\infty^n I(X_\infty>\theta)].\label{e.3.14} 
\end{equation}
Motivated from the approximations $\widehat L_{n,N}^\theta\approx L_n^\theta$ and  $\widehat R_{n,N}^\theta\approx R_n^\theta$,  we use  the following equations  to construct our   estimators for the parameters $\al_1, \al_2$: 
\begin{equation}
    \left\{
     \begin{aligned}
        &\frac{\widehat L_{1,N}^\theta}{\widehat L_{0,N}^\theta}  =-\frac{\sigma}{\sqrt{2\alpha_1}}\frac{\phi( -\sqrt{2\alpha_1}\theta/\sigma)}{1-\Phi(-\sqrt{2\alpha_1}\theta/\sigma)},\\
        &\frac{\widehat R_{1,N}^\theta}{\widehat R_{0,N}^\theta} =\frac{\sigma}{\sqrt{2\alpha_2}}\frac{\phi( \sqrt{2\alpha_2}\theta/\sigma)}{1-\Phi(\sqrt{2\alpha_2}\theta/\sigma)}.
     \end{aligned}
     \right.\label{e.3.15} 
\end{equation}
Let 
\begin{equation}
x=\frac{\sqrt{2\alpha_1}\theta}{\sigma}\,,\qquad y=\frac{\sqrt{2\alpha_2}\theta}{\sigma}\,, \label{e.3.16} 
\end{equation} 
and 
\begin{equation*}
A(x)=\frac{\phi(-x)}{1-\Phi(-x)}, \quad B(y)=\frac{\phi(y)}{1-\Phi(y)}. 
\end{equation*}
Equivalently, the system of equations \eqref{e.3.15} becomes   
\begin{equation}\label{e:L/L}
    \left\{
     \begin{aligned}
        &\frac{\widehat L_{1,N}^\theta}{\theta \widehat L_{0,N}^\theta} =-\frac{  A(x)}{x}=:K_1(x),\\
        &\frac{\widehat R_{1,N}^\theta}{\theta \widehat R_{0,N}^\theta}=\frac{  B(y)}{y}=:K_2(y).
     \end{aligned}
     \right.
\end{equation}
These are two uncoupled equations, so we can solve them separately.
To see if there is a unique solution to  each of the above equations or not, we use the simple mean value theorem: 
if a differentiable function $f$ has nonzero derivatives 
on an interval $I$,  then it is injective.  Using the fact that 
$A'(x)=-xA(x)-A^2(x)$ and $B'(y)=-yB(y)+B^2(y)$, 
we can compute the derivatives of $K_1$ and $K_2$ as follows: 
\begin{equation*}
    \left\{
     \begin{aligned}
        &\frac{dK_1}{d x}=   A(x)\left(\frac{1}{x^2}+1+\frac{A(x)}{x}\right),\\
        &\frac{dK_2}{d y}=-  B(y)\left (\frac{1}{y^2}+1-\frac{B(y)}{y}\right).
     \end{aligned}
     \right.
\end{equation*}
To investigate the monotonicity of $K_i$, $i=1,2$, it is equivalent to show the positivity or negativity of $F_1(x)=\frac{1}{x^2}+1+\frac{A}{x}$ and $F_2(y):=\frac{1}{y^2}+1-\frac{B}{y}$. Since $F_1(-y)=F_2(y)$, to show each of the 
equation in \eqref{e:L/L} has a unique solution in $\RR$,  we only need to show $F_1(x)>0$ for all $x\neq0$.  Denote $\widetilde F(x):=1-\Phi(-x)+x^2(1-\Phi(-x))+x\phi(-x)$. Then
  $F_1(x)=\widetilde F(x)/[x^2(1-\Phi(-x))]$.  Note that 
\begin{equation*}
\widetilde F^\prime(x)=2\phi(x)+2x\Phi(x),\quad \widetilde F^{\prime\prime}(x)=2\Phi(x)>0.
\end{equation*}
Since $\lim_{x\to-\infty}\widetilde F^\prime(x)=0$, we see  $\widetilde F^\prime(x)>0$. Now we can conclude that $\widetilde F(x)>0$ from  $\lim_{x\to-\infty}\widetilde F(x)=0$.
Therefore, there exists a continuous inverse function $H=(H_1, H_2)$ of $(K_1, K_2)$ such that
\begin{equation*}
\widehat x_N:=H_1\left(\frac{\widehat L_{1,N}^\theta}{ \th\widehat L_{0,N}^\theta} \right), \quad \widehat y_N:=H_2\left(\frac{\widehat R_{1,N}^\theta}{\th \widehat R_{0,N}^\theta} \right).
\end{equation*}
From the ergodic theorem   we know that  $\widehat L_{n,N}^\theta$ and $\widehat R_{n,N}^\theta$ converge almost surely to
$ L_{n }^\theta$ and $  R_{n }^\theta$ defined by \eqref{e.3.14}.    Thus,  the estimators $\widehat x_N$ and  $\widehat y_N $ converge almost surely to the parameters 
\begin{equation}\label{e:xy}
x=H_1(K_1(x))=\frac{\sqrt{2\alpha_1}\theta}{\sigma}, \quad  y=H_2(K_2(y))=\frac{\sqrt{2\alpha_2}\theta}{\sigma}
\end{equation}  
respectively, as $N\to\infty$. 
Now the relationship \eqref{e.3.16} 
 between $(x,y)$ and $(\al_1,\al_2)$  yields the following theorem. 
\begin{thm}\label{thm:3.3} 
For any sample size $N$ the system of equations \eqref{e:L/L} has a unique solution $(\widehat x_N,\widehat y_N)$. 
The generalized moment estimators defined by
\begin{equation*}
 \widehat \alpha_{1,N}=\frac{1}{2}\left(\frac{\sigma\widehat x_N}{\theta} \right)^2, \quad 
 \widehat \alpha_{2,N}=\frac{1}{2}\left(\frac{\sigma\widehat y_N}{\theta} \right)^2 
\end{equation*}
 are strongly consistent, namely,   
 $(\widehat \alpha_{1,N} ,   
 \widehat \alpha_{2,N})$ converges to $(\alpha_1,\alpha_2)$
 almost surely.
\end{thm}
Compared with the case I, the estimators only have implicit expressions in terms of the inverse functions $H_1$ and $H_2$. Nevertheless, it is clear that $H_1$ and $H_2$ are continuously differentiable. Hence, we can exhibit the following  CLT for the estimators $\widehat\alpha_{i,N}$, $i=1,2$.
\begin{thm}\label{t.3.7} 
As $N\to\infty$,
\begin{equation*}
\sqrt N \left((\widehat \alpha_{1,N},\widehat \alpha_{2,N})^T-(\alpha_1,\alpha_2)^T\right)\Rightarrow \mathbf N(0,\widehat\Sigma),
\end{equation*}
where $\widehat\Sigma$ is given by \eqref{e.3.22}  below. 
\end{thm}
\begin{pf}
The proof is similar to that of Theorem \ref{asy:1}, so we only provide a sketch of the proof.  
Set
\begin{align*}
&F_1(x)=I(x\le\theta), \quad F_2(x)=xI(x\le\theta) \,, \\
&F_3(x)=I(x>\theta), \quad F_4(x)=xI(x>\theta) \,. 
\end{align*}
From \citet[Theorem 17.0.1]{meyn2012markov}, we get that for $i,j=1,2,3,4,$
\begin{equation*}
\widetilde\sigma_{ij}:=\Cov (F_{i}(\widetilde X_0),F_{j}(\widetilde X_0))+\sum_{k=1}^\infty\left[ \Cov (F_{i}(\widetilde X_0),F_{j}(
\widetilde X_{kh}))+\Cov (F_{j}(\widetilde X_0),F_{i}(\widetilde X_{kh})) \right],
\end{equation*}
are well defined and  non-negative.  They can be computed by using \eqref{e.B.2} as follows:
\begin{equation*}
\tilde \si_{ij}=\sigma(F_i, F_j)\,, \quad i,j=1,  2, 3, 4\,. 
\end{equation*}
Denote $\tilde\Sigma_2:=(\tilde\sigma _{ij})_{1\le i,j\le 4}$,
 then we have 
\begin{equation*}
\sqrt N\left(( \widehat L_{0,N}^\theta, \widehat L_{1,N}^\theta, \widehat R_{0,N}^\theta, \widehat R_{1,N}^\theta)^T-(  L_{0}^\theta, L_{1}^\theta, R_{0}^\theta, R_{1}^\theta) ^T\right)\Rightarrow \mathbf N(0,\tilde\Sigma_2).
\end{equation*}
Define two functions by $h_1(x_1,x_2):= H_1(\frac{x_2}{\theta x_1} )$ and $h_2(x_3,x_4):=H_2(\frac{x_4}{\theta x_3})$ and set two maps
\begin{align*}
&h:(x_1,x_2,x_3,x_4)\mapsto  (h_1(x_1,x_2),h_2(x_3,x_4)) \;, \\
& l:(x_1,x_2) \mapsto \left(\frac{\sigma^2 {x_1}^2}{2\theta^2}, \frac{\sigma^2 {x_2}^2}{2\theta^2}\right)\; .
\end{align*}
By the multivariate delta method, we have
\begin{equation*}
\sqrt N\left(h(\widehat L_{0,N}^\theta, \widehat L_{1,N}^\theta, \widehat R_{0,N}^\theta, \widehat R_{1,N}^\theta)^T-h( L_{0}^\theta, L_{1}^\theta, R_{0}^\theta, R_{1}^\theta)^T\right) \Rightarrow \mathbf N(0, \bar\Sigma ),
\end{equation*}
where $\bar\Sigma=\nabla h( L_{0}^\theta, L_{1}^\theta, R_{0}^\theta, R_{1}^\theta)\tilde\Sigma_2  \nabla h( L_{0}^\theta, L_{1}^\theta, R_{0}^\theta, R_{1}^\theta)^T $.  Applying the multivariate delta method again, we get the desired CLT result
\begin{equation*}
\sqrt N\left(l(h(\widehat L_{0,N}^\theta, \widehat L_{1,N}^\theta, \widehat R_{0,N}^\theta, \widehat R_{1,N}^\theta))^T-l(h( L_{0}^\theta, L_{1}^\theta, R_{0}^\theta, R_{1}^\theta))^T\right) \Rightarrow \mathbf N(0, \widehat\Sigma ),
\end{equation*}
where 
\begin{equation}
\widehat\Sigma:=\nabla l(h( L_{0}^\theta, L_{1}^\theta, R_{0}^\theta, R_{1}^\theta) )\ \bar\Sigma \ \nabla l(h( L_{0}^\theta, L_{1}^\theta, R_{0}^\theta, R_{1}^\theta) )^T\,.  \label{e.3.22}
\end{equation} 
 The proof
is then completed.
\end{pf}


\subsection{Case III: Estimate $\beta_i$ and $\alpha_i$ for known $\theta\neq0$}\label{sub:III}
\begin{figure} 
\centering 
\subfigure[$D_1(u)$]{\label{fig:a1}
\includegraphics[width=0.6\textwidth]{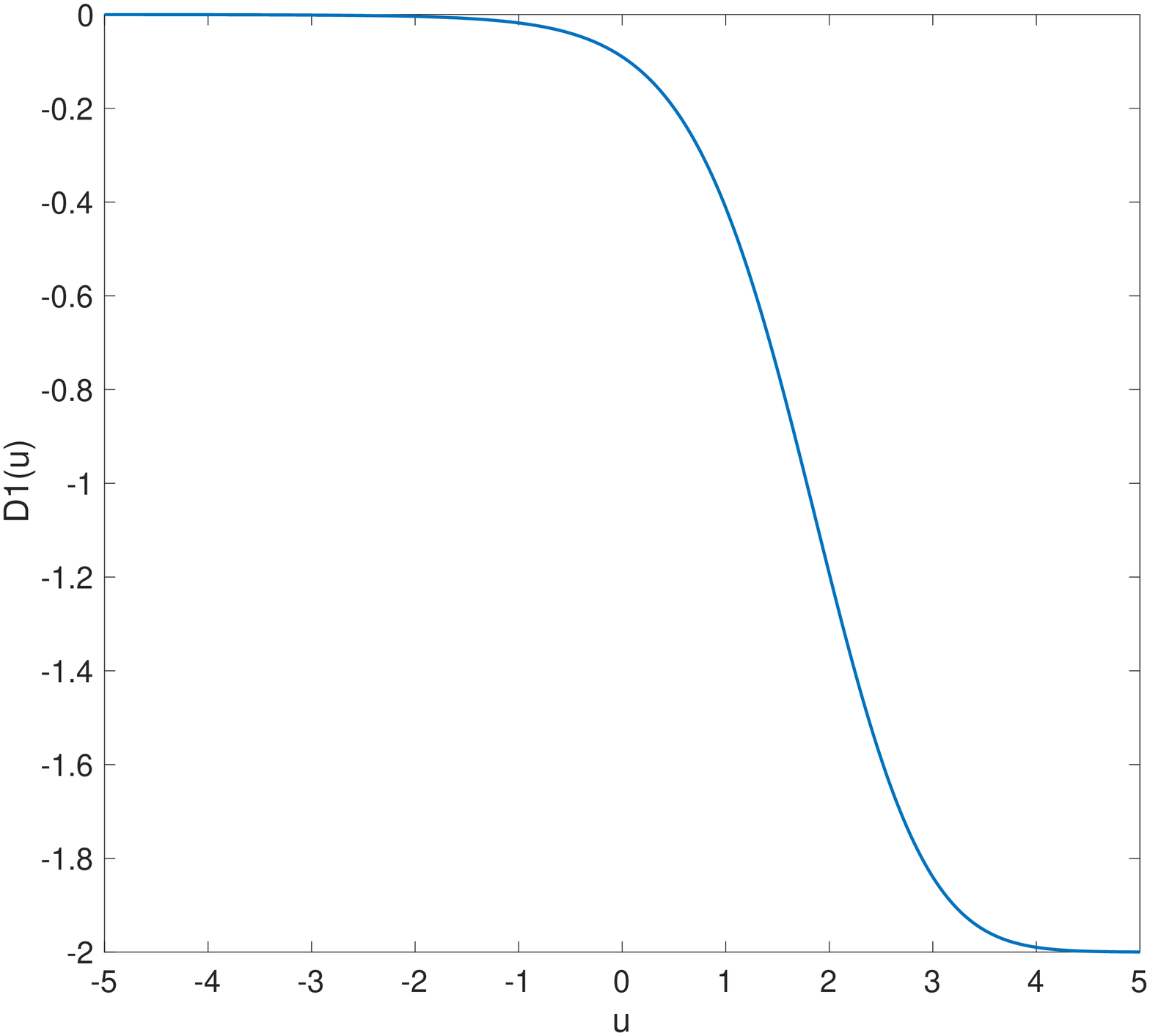}}
\subfigure[ $D_1(u)$]{\label{fig:b1}
\includegraphics[width=0.6\textwidth]{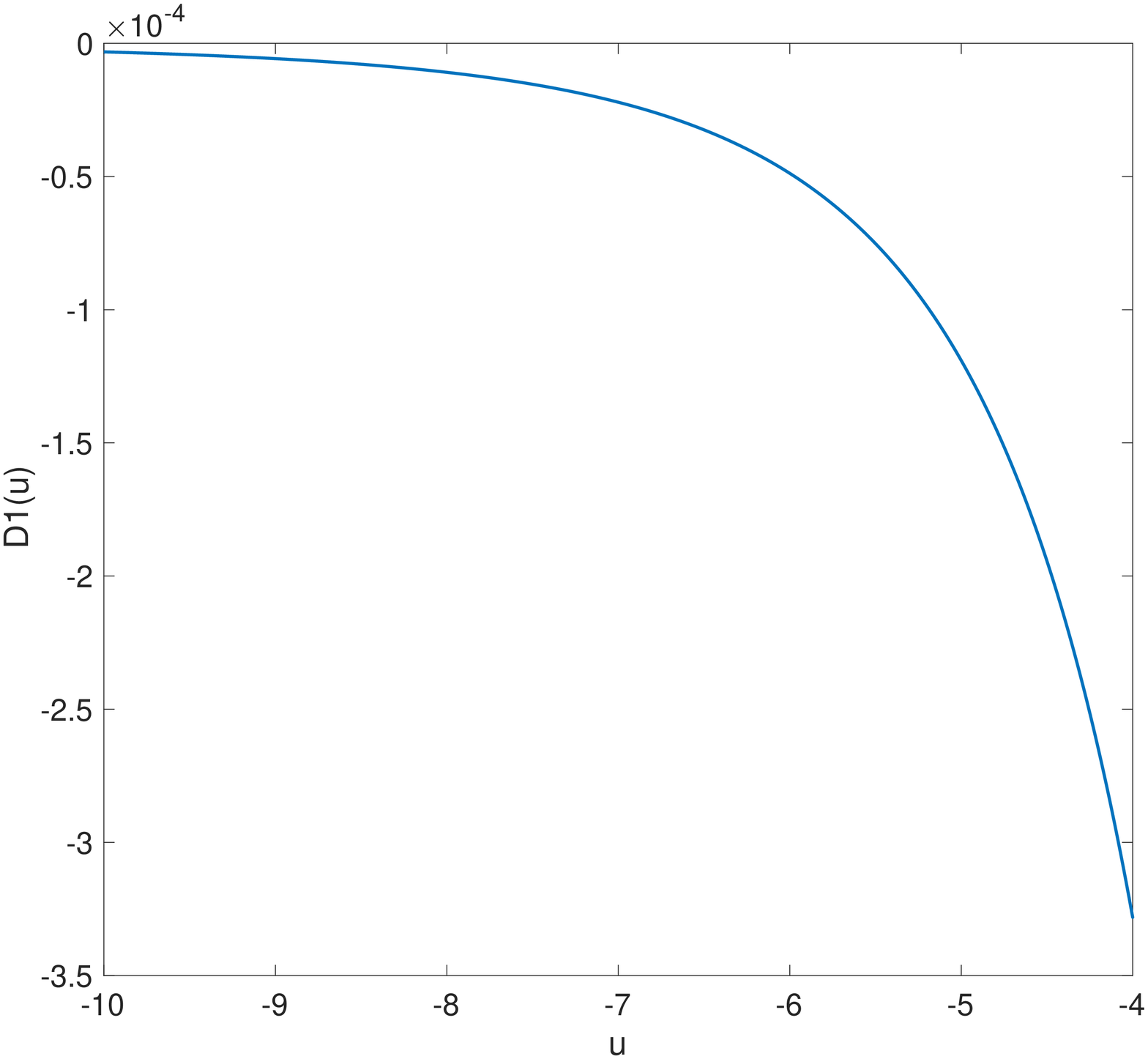}}
\caption{The plot of $D_1(u)$. }
\label{fig:D1}
\end{figure}

In this subsection,  we extend our approach to multiple-parameter case, where $\theta\neq 0$. The stationary density is given by
\begin{equation}
\psi_3(x)=k_1\exp\left(\frac{-\alpha_1x^2+2\beta_1x}{\sigma^2}\right)I(x\le \theta)+k_2\exp\left(\frac{-\alpha_2x^2+2\beta_2x}{\sigma^2}\right)I(x> \theta),\label{e.general_psi}
\end{equation}
where $k_1$ and $k_2$ are defined by \eqref{e:k1} and    
\eqref{e:k2}. 
We can obtain the following stationary moments
\begin{equation}
    \left\{
     \begin{aligned}
        &\frac{\EE [X_\infty I(X_\infty\le\theta)]}{\EE [I(X_\infty\le\theta)]}=\frac{-\frac{\sigma}{\sqrt{2\alpha_1}}\phi\left(\frac{\sqrt{2\alpha_1}\theta}{\sigma}-\frac{2\beta_1}{\sqrt{2\alpha_1}\sigma}\right)}{\Phi\left(\frac{\sqrt{2\alpha_1}\theta}{\sigma}-\frac{2\beta_1}{\sqrt{2\alpha_1}\sigma}\right)}+\frac{\beta_1}{\alpha_1},\\
        &\frac{\EE [ X_\infty ^2 I(X_\infty\le\theta)]}{\EE [I(X_\infty\le\theta)]} =\frac{\sigma^2}{2\alpha_1}+\left(\frac{\beta_1}{\alpha_1}\right)^2 +\frac{\phi\left(\frac{\sqrt{2\alpha_1}\theta}{\sigma}-\frac{2\beta_1}{\sqrt{2\alpha_1}\sigma}\right)}{ \Phi\left(\frac{\sqrt{2\alpha_1}\theta}{\sigma}-\frac{2\beta_1}{\sqrt{2\alpha_1}\sigma}\right)} \left( -\theta-\frac{\beta_1}{\alpha_1}\right) \frac{\sigma}{\sqrt{2\alpha_1}}\,, \\
        &\frac{\EE [X_\infty I(X_\infty>\theta)]}{\EE [I(X_\infty>\theta)]}=\frac{\frac{\sigma}{\sqrt{2\alpha_2}}\phi\left(\frac{\sqrt{2\alpha_2}\theta}{\sigma}-\frac{2\beta_2}{\sqrt{2\alpha_2}\sigma}\right)}{1-\Phi\left(\frac{\sqrt{2\alpha_2}\theta}{\sigma}-\frac{2\beta_2}{\sqrt{2\alpha_2}\sigma}\right)}+\frac{\beta_2}{\alpha_2},\\
        &\frac{\EE [ X_\infty ^2 I(X_\infty>\theta)]}{\EE [I(X_\infty>\theta)]}=\frac{\sigma^2}{2\alpha_2}+\left(\frac{\beta_2}{\alpha_2}\right)^2 +\frac{\phi\left(\frac{\sqrt{2\alpha_1}\theta}{\sigma}-\frac{2\beta_1}{\sqrt{2\alpha_1}\sigma}\right)}{1- \Phi\left(\frac{\sqrt{2\alpha_1}\theta}{\sigma}-\frac{2\beta_1}{\sqrt{2\alpha_1}\sigma}\right)} \left( \theta+\frac{\beta_1}{\alpha_1}\right) \frac{\sigma}{\sqrt{2\alpha_2}}\,. 
     \end{aligned}
     \right. \label{e.3.20} 
\end{equation}
Denote the right-hand sides of the above identities by $\bar K_i$, $i=1,2,3,4$. Let 
\begin{align}\label{u,v}
v&=\frac{\beta_1}{\alpha_1}, \quad u=\frac{\sqrt{2\alpha_1}\theta}{\sigma}-\frac{2\beta_1}{\sqrt{2\alpha_1}\sigma}=\frac{\sqrt{2\alpha_1}(\theta-v)}{\sigma},\quad  A(u)=\frac{\phi(-u)}{1-\Phi(-u)}, \\\label{w,z}
z&=\frac{\beta_2}{\alpha_2}, \quad \omega=\frac{\sqrt{2\alpha_2}\theta}{\sigma}-\frac{2\beta_2}{\sqrt{2\alpha_2}\sigma}=\frac{\sqrt{2\alpha_2}(\theta-z)}{\sigma}, \quad B(\om)=\frac{\phi(\omega)}{1-\Phi(\omega)}\,. 
\end{align}
Then we can rewrite $\bar K_i$ as
\begin{equation}
    \left\{
     \begin{aligned}
        &\bar K_1(u,v) =\frac{v-\theta}{u}A(u)+v,\\
        &\bar K_2(u,v)=\left( \frac{\theta-v}{u}\right)^2+v^2-A(u) \frac{\theta^2-v^2}{u}\,,  \\
        &\bar K_3(\omega,z)=\frac{\theta-z}{\omega}B(\om)+z,\\
        &\bar K_4(\omega,z)=\left( \frac{\theta-z}{\omega}\right)^2+z^2+B(\om)\frac{\theta^2-z^2}{\omega}\,. 
     \end{aligned}
     \right. \label{e.3.23} 
\end{equation}
 Similar to  the previous cases,  we approximate the left hand sides of \eqref{e.3.20} by the following statistics for $i=1,2$: 
 \begin{align*}
 \widehat L_{i,N}^\theta/\widehat L_{0,N}^\theta&\approx \EE[(X_\infty)^i I(X_\infty\le\theta)]/ \EE[ I(X_\infty\le\theta)],\\
 \widehat R_{i,N}^\theta/\widehat R_{0,N}^\theta&\approx \EE[(X_\infty)^i I(X_\infty\le\theta)]/ \EE[ I(X_\infty>\theta)], 
 \end{align*}
Motivated by \eqref{e.3.20}   and \eqref{e.3.23} 
  we first propose the 
following    estimators $\widehat v_N$, $\widehat u_N$, $\widehat z_N$, and $\widehat \omega_N$  to estimate  $v, u, z, \om$ 
by solving the following system
\begin{equation}\label{sys:III}
    \left\{
     \begin{aligned}
        &\frac{\widehat L_{1,N}^\theta}{\widehat L_{0,N}^\theta}  =\frac{v-\theta}{u}A(u)+v,\\
        &\frac{\widehat L_{2,N}^\theta}{\widehat L_{0,N}^\theta} =\left( \frac{\theta-v}{u}\right)^2+v^2-A(u)\frac{\theta^2-v^2}{u}\,,  \\
        &\frac{\widehat R_{1,N}^\theta}{\widehat R_{0,N}^\theta} =\frac{\theta-z}{\omega}B(\om)+z,\\
        &\frac{\widehat R_{2,N}^\theta}{\widehat R_{0,N}^\theta} =\left( \frac{\theta-z}{\omega}\right)^2+z^2+B(\om)\frac{\theta^2-z^2}{\omega}\,. 
     \end{aligned}
     \right.
\end{equation}
Next we need to solve this system   of four equations. 
First,  we observe  that 
this  system of four  equations is decoupled as two systems, each consisting two equations. Let us first study the
first pair of  equations in \eqref{sys:III}: 
\begin{equation}\label{sys:III_a}
    \left\{
     \begin{aligned}
        &\frac{\widehat L_{1,N}^\theta}{\widehat L_{0,N}^\theta}  =\frac{v-\theta}{u}A(u)+v=:\bar K_1(u,v)\,, \\
        &\frac{\widehat L_{2,N}^\theta}{\widehat L_{0,N}^\theta} =\left( \frac{\theta-v}{u}\right)^2+v^2-A(u)\frac{\theta^2-v^2}{u}=:\bar K_2(u,v)\,. 
     \end{aligned}
     \right.
\end{equation}
 The partial  derivatives of $\bar K_1, \bar K_2$ are given by
\begin{equation*}
    \left\{
     \begin{aligned}
        &\frac{\partial \bar K_1}{\partial u} = -\frac{v-\theta}{u^2}A(u)-(v-\theta)A(u)-A^2(u)\frac{v-\theta}{u} \,, \\
        &\frac{\partial \bar K_1}{\partial v}= \frac{A(u)}{u}+1\,,\\
        &\frac{\partial \bar K_2}{\partial u}=-\frac{2(\theta-v)^2}{u^3}- (-uA(u)-A^2(u))\frac{\theta^2-v^2}{u}+A(u)  \frac{\theta^2-v^2}{u^2}\,,\\
        &\frac{\partial \bar K_2}{\partial v}=-\frac{2(\theta-v)}{u^2}+2v+\frac{2A(u)v}{u}\,. 
     \end{aligned}
     \right.
\end{equation*}
The Jacobian matrix $J_1$ of $(\bar K_1, \bar K_2)$ is given by
\begin{equation*}
J_1= \begin{pmatrix}
\frac{\partial \bar K_1}{\partial u}&\frac{\partial \bar K_1}{\partial v}\\
\frac{\partial \bar K_2}{\partial u}&\frac{\partial \bar K_2}{\partial v}\\
\end{pmatrix}.
\end{equation*}  
The determinant of $J_1$ is 
\begin{align*}
\det(J_1)=-\frac{(v-\theta)^2}{u^3}(A(u)u^3+3A(u)u+A^3(u)u+2A^2(u) u^2+3A^2(u)-2) \;.
\end{align*}
  Let $D_1(u)=A(u)u^3+3A(u)u+A^3(u)u+2A^2(u)u^2+3A^2(u)-2$. To show that $\det(J_1)\neq0$ for all $u\neq0$ and  $v\neq\theta$ it suffices to show that $D_1(u)<0$.
  From  the Figure \ref{fig:D1}, we can see that $D_1(u)<0$  for all $u\in[-10,5]$.
Let 
\begin{eqnarray*}
\DD_1&=&\left\{(u,v)\in \RR^2\,; v\not =\theta\,\right.
 \   {\rm and}  \nonumber\\
&&\left. D_1(u)=A(u)u^3+3A(u)u+A^3(u)u+2A^2(u)u^2+3A^2(u)-2 \not =0 \right\}\,.  
\end{eqnarray*} 
The Figure  \ref{fig:D1} implies $  \left\{
(u, v), u\in (-10, 5)\,, v\not =0\right\}\subseteq \DD_1$. 
If necessary, one can enlarge the interval $(-10, 5)$.   
If  $(u_0,v_0)\in \DD_1$ is from the   true parameters 
$(\al_1, \be_1)$, then by the ergodic Lemma \ref{le:ergodic}
we know  that when  $N$ goes to infinity  \eqref{sys:III_a}
will become  true identities with the $(u,v)$ on the right-hand side replaced by $(u_0, v_0)$. 
Thus, when $N$ is sufficiently large
$\displaystyle \left(\frac{\widehat L_{1,N}^\theta}{\widehat L_{0,N}^\theta}\,, \frac{\widehat L_{2,N}^\theta}{\widehat L_{0,N}^\theta} \right)$ will be in any given neighbourhood of 
$(\bar K_1(u_0,v_0), \bar K_2(u_0,v_0))$. 
On the other hand, it is obvious that $\DD_1$ is an open set in $\RR^2$ and 
$\bar K_1, \bar K_2$ are continuous functions of $(u,v)$. 
Since  
$\det(J_1)\not=0$ on $\DD_1$, by the inverse function theorem
there is one unique solution pair $(u,v)\in \DD_1$ in some neighbourhood of $(u_0, v_0)$ such that the system of
equations \eqref{sys:III_a}  are satisfied. 
This gives the existence and local uniqueness of the 
solution to the system of
equations \eqref{sys:III_a}. 

\begin{figure}[t] 
\centering 
\includegraphics[width=0.6\textwidth]{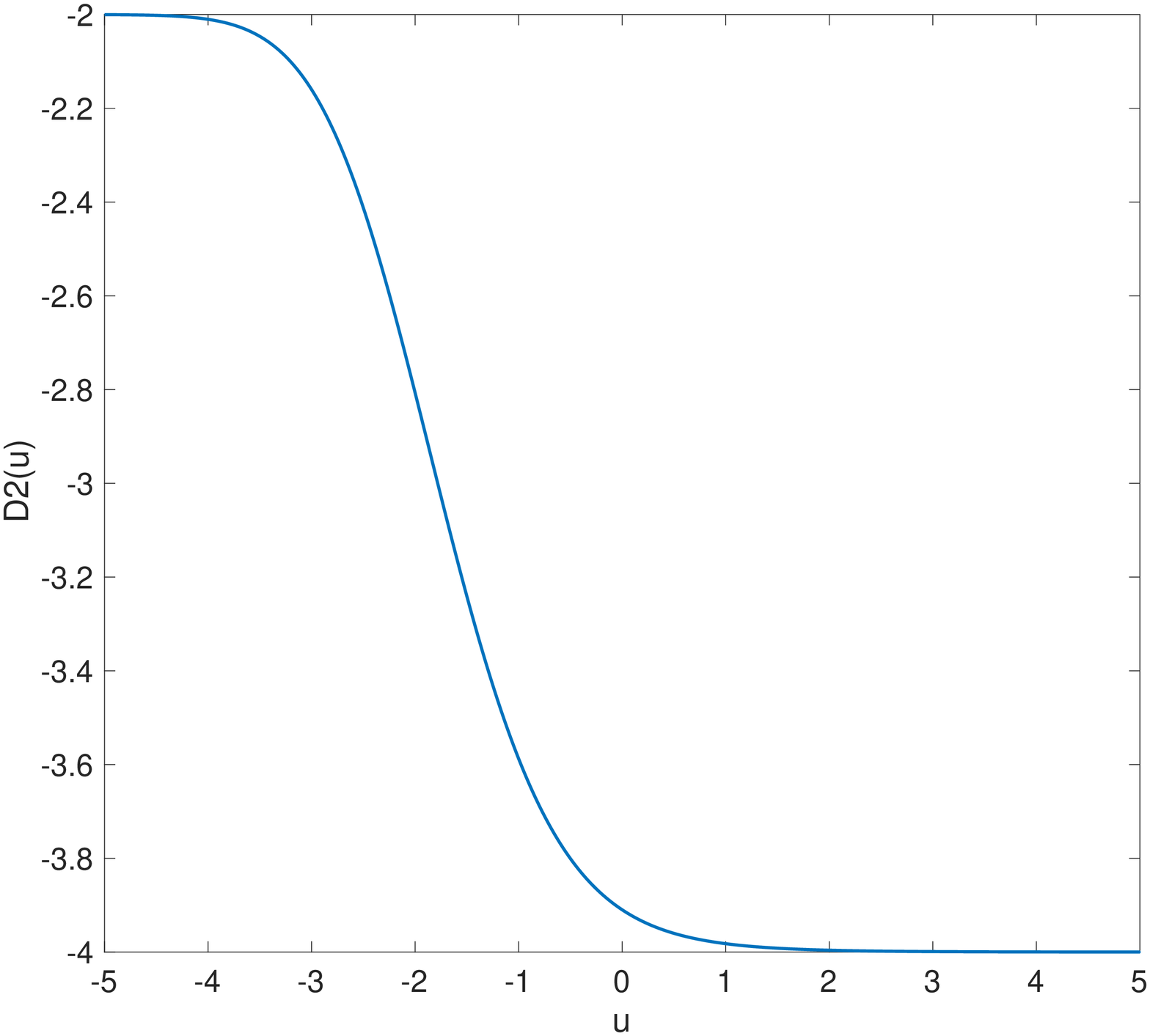} 
\caption{$D_2(u)$} 
\label{Fig.D2} 
\end{figure}

Now we consider the second pair of   equations in \eqref{sys:III}. 
\begin{equation}\label{sys:III_b}
    \left\{
     \begin{aligned}
        &\frac{\widehat R_{1,N}^\theta}{\widehat R_{0,N}^\theta} =\frac{\theta-z}{\omega}B(\om)+z=:\bar K_3(  \om, z)\,,\\
        &\frac{\widehat R_{2,N}^\theta}{\widehat R_{0,N}^\theta} =\left( \frac{\theta-z}{\omega}\right)^2+z^2+B(\om)\frac{\theta^2-z^2}{\omega}=:\bar K_4(  \om, z) \,. 
     \end{aligned}
     \right.
\end{equation}
The partial derivatives of $\bar K_3(\om, z)$, $\bar  K_4(\om, z)$ are 
\begin{equation*}
    \left\{
     \begin{aligned}
        &\frac{\partial \bar K_3}{\partial  \omega} = -\frac{\theta-z}{\omega^2}B(\om)-(\theta-z)B(\om)+B^2(\om)\frac{\theta-z}{\omega}\\
        &\frac{\partial \bar K_3}{\partial z}=-\frac{B(\om)}{\omega}+1,\\
        &\frac{\partial \bar K_4}{\partial \omega}=-\frac{2(\theta-z)^2}{\omega^3}+(-\omega B(\om)+B^2(\om))\frac{\theta^2-z^2}{\omega}-B(\om) \frac{\theta^2-z^2}{\omega^2},\\
        &\frac{\partial \bar K_4}{\partial z}=-\frac{2(\theta-z)}{\omega^2}+2z-\frac{2B(\om)z}{\omega}\\
     \end{aligned}
     \right.
\end{equation*}
The determinant of the Jacobian matrix $J_2$ of $(\bar K_3, \bar K_4)$  is 
\begin{align*}
\det(J_2)=\frac{(\theta-z)^2}{\omega^3}(B(\om)\omega^3+3B(\om)\omega+B^3(\om)\omega-2B^2(\om)\omega^2-3B^2(\om)-2)\;.
\end{align*}
Let $D_2(\omega)=B(\om)\omega^3+3B(\om)\omega+B^3(\om)\omega-2B^2(\om)\omega^2-3B^2(\om)-2$.  From the Figure \ref{Fig.D2},  we see that $D_2(\omega)<0$ for all $\om\in [-5, 5]$. 
Denote 
\begin{eqnarray*}
\DD_2&=&\left\{(z, \om )\in \RR^2\,; z\not =\theta\,\right.
 \   {\rm and}  \nonumber\\
&&\left. D_2(\omega)=B(\om)\omega^3+3B(\om)\omega+B^3(\om)\omega-2B^2(\om)\omega^2-3B^2(\om)-2 \not =0 \right\}\,.  
\end{eqnarray*} 
Analogous to the argument for the system of
equations \eqref{sys:III_a} we can prove  the existence and local uniqueness of the 
solution to the system of
equations \eqref{sys:III_b}. 

Once we have the existence and local uniqueness of
 the system of
equations \eqref{sys:III_a} and \eqref{sys:III_b} we can 
follow    the substitutions \eqref{u,v} and \eqref{w,z} to  obtain the generalized moment estimators  $\widehat\alpha_{i,N}$ and $\widehat\beta_{i,N}$  for $\al_i$ and $\be_i$, $i=1, 2$.
We summarize the above as the following theorem.
\begin{thm}\label{t.3.8} 
Let $(\al_1, \be_1, \al_2, \be_2)$ be the true parameters such that $(u,v)$ and $(z, \om)$  defined  by \eqref{u,v} and \eqref{w,z}   are in $\DD_1$ and $\DD_2$, respectively.  Then,  when $N$ is sufficiently large 
the systems of equations  \eqref{sys:III_a} and \eqref{sys:III_b} have   solutions  $(\widehat u_N, \widehat v_N)$
and $( \widehat \om_N,\widehat z_N)$, respectively. The solutions are unique in a neighbourhood of $(u,v)$ and a neighbourhood of
$(\om, z)$. If we define 
\begin{empheq}[left=\empheqlbrace]{align} 
\widehat \alpha_{1,N}&=\frac{(\widehat u_N)^2\sigma^2}{2(\theta-\widehat v_N)^2}, \quad \widehat \alpha_{2,N}=\frac{(\widehat \omega_N)^2\sigma^2}{2(\theta-\widehat z_N)^2},\label{e.3.33}\\
\widehat \beta_{1,N}&=\widehat v_N \widehat\alpha_{1,N}, \quad \widehat \beta_{2,N}=\widehat z_N \widehat \alpha_{2,N}\,,  
\label{e.3.34}
\end{empheq}
then when $N\rightarrow \infty$,  we have
\begin{equation*}
(\widehat \alpha_{1,N}\,, \widehat \alpha_{2,N}\,, 
\widehat \beta_{1,N}\,,  \widehat \beta_{2,N} )\rightarrow (\al_1, \al_2, \be_1, \be_2)
\quad \hbox{almost suely}\,. 
\end{equation*}
\end{thm}
%
\begin{rmk}\label{rmk:3.10}
If $u=0$, $\omega=0$, i.e., $\frac{\beta_i}{\alpha_i}=v=z=\theta$, $i=1,2$. We can estimate $(\alpha_1,\alpha_2)$ by solving
\begin{equation*}
    \left\{
     \begin{aligned}
        &\frac{\widehat L_{1,N}^\theta}{\widehat L_{0,N}^\theta} = -\frac{\sigma}{\sqrt{\pi\alpha_1}}+\theta,\\
        &\frac{\widehat R_{1,N}^\theta}{\widehat R_{0,N}^\theta} = \frac{\sigma}{\sqrt{\pi\alpha_2}}+\theta.\\
     \end{aligned}
     \right.
\end{equation*}
Then $\widehat\beta_i=\theta\widehat\alpha_i$, $i=1,2$.
\end{rmk}
We also have the CLT for the above estimators.
Before stating the theorem, let us describe the   asymptotic variances.  Let 
\begin{empheq}[left=\empheqlbrace]{align*} 
&G_1(x)=I(x\le\theta), \quad G_2(x)=xI(x\le\theta) , \quad G_3(x)=x^2I(x\le\theta)\,, \\
&G_4(x)=I(x>\theta), \quad G_5(x)=xI(x>\theta)
, \quad G_6(x)=x^2I(x>\theta) \,.
\end{empheq}
Denote 
\begin{equation*}
\tilde\Sigma_3=\left(\si_{ij}\right)_{1\le i,j\le 6}\,,
\quad \hbox{where}\quad \si_{ij}=\sigma(G_i, G_j)\,, 1\le i,j\le 6 
\end{equation*}
with $\sigma(G_i, G_j)$ being defined by \eqref{e.B.2}. 
Then we have as before, 
\begin{equation*}
\sqrt N\left( ( \widehat L_{0,N}^\theta, \widehat L_{1,N}^\theta,\widehat L_{2,N}^\theta, \widehat R_{0,N}^\theta, \widehat R_{1,N}^\theta, R_{2,N}^\theta)^T-(  L_{0}^\theta, L_{1}^\theta,L_{2}^\theta, R_{0}^\theta, R_{1}^\theta,  R_{2}^\theta) ^T\right)\Rightarrow \mathbf N(0, \tilde\Sigma_3).
\end{equation*}
Let $(u,v)  =(\kappa_1 (x_1, x_2), \kappa_2(x_1, x_2))$ be the inverse mapping of $(\bar K_1(u,v), \bar K_2(u,v))$  defined by \eqref{sys:III_a} and 
let $(\om,z )  =(\kappa_3 (x_3, x_4), \kappa_4(x_3, x_4))$ be the inverse mapping of $(\bar K_3(u,v), \bar K_4(u,v))$  defined by \eqref{sys:III_b}.  Comparing with \eqref{e.3.33}-\eqref{e.3.34}  and denoting $x=(x_1, x_2, x_3, x_4, x_5, x_6)$,  we introduce 
\begin{empheq}[left=\empheqlbrace]{align} 
&\rho_1(x):= \frac{(\kappa_1 (\frac{x_2}{x_1}, \frac{x_3}{x_1}))^2\sigma^2}{2(\theta-\kappa_2 (\frac{x_2}{x_1}, \frac{x_3}{x_1}))^2}\,;\nonumber\\
 &\rho_2(x):=\frac{(\kappa_3(\frac{x_5}{x_4}, \frac{x_6}{x_4})^2\sigma^2}{2(\theta-\kappa_4(\frac{x_5}{x_4}, \frac{x_6}{x_4}))^2}\,;\nonumber\\ 
&\rho_3(x):=\kappa_2(\frac{x_2}{x_1}, \frac{x_3}{x_1}) \rho_1(x)\,;\nonumber\\
&\rho_4(x):=\kappa_4(\frac{x_5}{x_4}, \frac{x_6}{x_4})  \rho_2(x) \,.   \nonumber
\end{empheq}
Define a map $\rho:\RR^6\ni x\mapsto(\rho_1(x),\rho_2(x),\rho_3(x),\rho_4(x))\in \RR^4$. Now we establish the following asymptotic normality theorem.
\begin{thm}\label{t.3.10} 
As $N\to \infty$, we have the following asymptotic normality:
\begin{equation*}
\sqrt N\left( (\widehat \alpha_{1,N}\,, \widehat \alpha_{2,N}\,, 
\widehat \beta_{1,N}\,,  \widehat \beta_{2,N} )^T - (\al_1, \al_2, \be_1, \be_2)^T\right)\Rightarrow \mathbf
N(0, \bar  \Sigma_3  )\,, 
\end{equation*}
where 
\begin{equation*}
\bar  \Sigma_3 =\nabla \rho
(  L_{0}^\theta, L_{1}^\theta,L_{2}^\theta, R_{0}^\theta, R_{1}^\theta,  R_{2}^\theta)\ \tilde\Sigma_3  \nabla \rho
(  L_{0}^\theta, L_{1}^\theta,L_{2}^\theta, R_{0}^\theta, R_{1}^\theta,  R_{2}^\theta)^T\,. 
\end{equation*} 
\end{thm} 
Theorem \ref{t.3.8} gives domains $\DD_1$ and $\DD_2$ so that
we can find generalized moment estimators
$\widehat \alpha_{1,N}\,, \widehat \alpha_{2,N}$\,,  $
\widehat \beta_{1,N}\,,  \widehat \beta_{2,N}$ of $ \al_1, \al_2, \be_1, \be_2 $.  
On the one hand, although the functions $D_1$ and $D_2$ are 
explicit, we still have difficulty to know the shapes of  $\DD_1$ and $\DD_2$. Our numerical experiments suggest 
that $D_1(u) \not =0$ and $D_2(u)\not =0$ for all $u\in 
\RR$.  However,  we cannot conclude this analytically. 
On the other hand, as we know that the implicit function
theorem is a local one in high dimensions. 
This means that the solutions to \eqref{sys:III_a} and
to  \eqref{sys:III_b}  are unique only in a neighbourhood of
the true parameters. The method of nondegeneracy  of the 
determinant cannot be used to guarantee the existence of a global inverse function. For example, the  mapping 
$(f(x,y), g(x,y))=(e^x\cos y, e^x \sin y)$  from $\RR^2$ to $\RR^2$ has   a  strictly 
positive 
Jacobian determinant $J(f, g)=e^x$ on the whole plane $\RR^2$. But it is not an injection as a mapping from $\RR^2$ to $\RR^2$.  
Therefore,  Theorem \ref{t.3.8} is powerful when we know a priori roughly the range of the true parameters.  For example,
in the modelling of the financial market, we know roughly the 
long memory Hurst parameter $H$ is around $0.5$. But in some other cases  researchers do not have any idea about the parameter ranges. 
Thus, a natural question arises: What should we do if there are more than one solution to \eqref{sys:III_a} and
to  \eqref{sys:III_b}? 
Now we are going to address this global uniqueness issue
(existence is not an issue by Theorem \ref{t.3.8}).

From the first equation of \eqref{sys:III_a} we have
\begin{equation}
v=\frac{u\frac{\widehat L_{1,N}^\theta}{\widehat L_{0,N}^\theta}+\theta A(u)}{u+A(u)}\,. \label{e.3.32}
\end{equation}
Substituting it to the second equation of \eqref{sys:III_a} we  obtain
\begin{eqnarray}
&&u\left(\frac{\widehat L_{1,N}^\theta}{\widehat L_{0,N}^\theta}
 -\theta\right)^2+ u\left(u\frac{\widehat L_{1,N}^\theta}{\widehat L_{0,N}^\theta}
+\theta A(u)\right)^2\nonumber\\
&&\qquad\qquad-A(u) \left[\theta^2(u+A(u))^2-
\left(u\frac{\widehat L_{1,N}^\theta}{\widehat L_{0,N}^\theta}
+\theta A(u)\right)^2\right] - u(u+A(u))^2
\frac{\widehat L_{2,N}^\theta}{\widehat L_{0,N}^\theta} =0\,. 
\end{eqnarray} 
This is one  equation on one unknown $u$.  
Solving $F_1(u)=0$ to get $\widehat u_N$ and substituting 
it into \eqref{e.3.32}, we can get $\widehat  v_N$.
Notice that the 
quantities $\frac{\widehat L_{1,N}^\theta}{\widehat L_{0,N}^\theta}$ and $\frac{\widehat L_{2,N}^\theta}{\widehat L_{0,N}^\theta}$ appeared in \eqref{e.3.33} can be computed from real data. 

\begin{figure} [t]
\centering 
\subfigure[$F_1(u)$]{\label{fig:F1}
\includegraphics[width=0.48\textwidth]{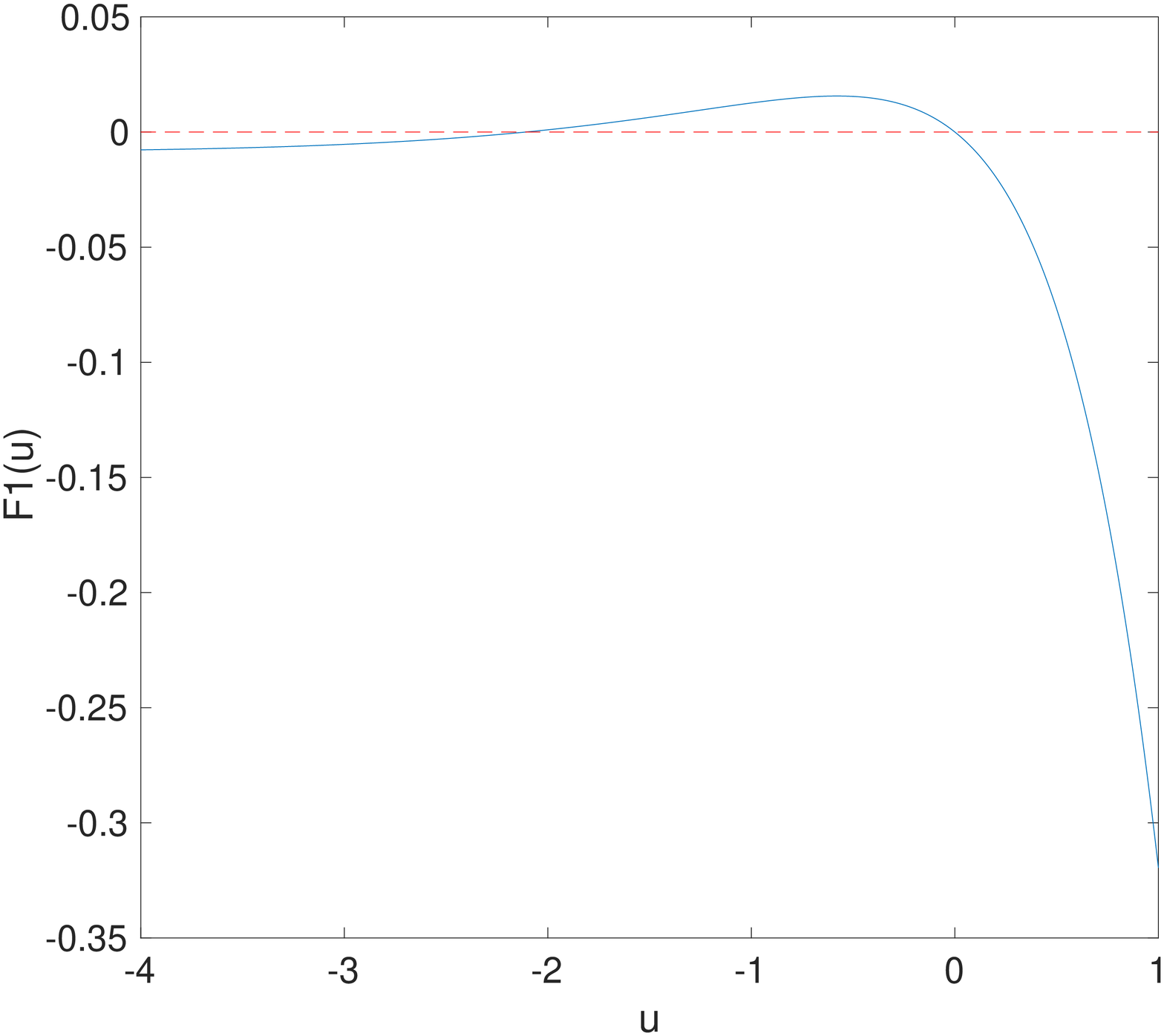}}
\subfigure[ $F_2(\om)$]{\label{fig:F2}
\includegraphics[width=0.48\textwidth]{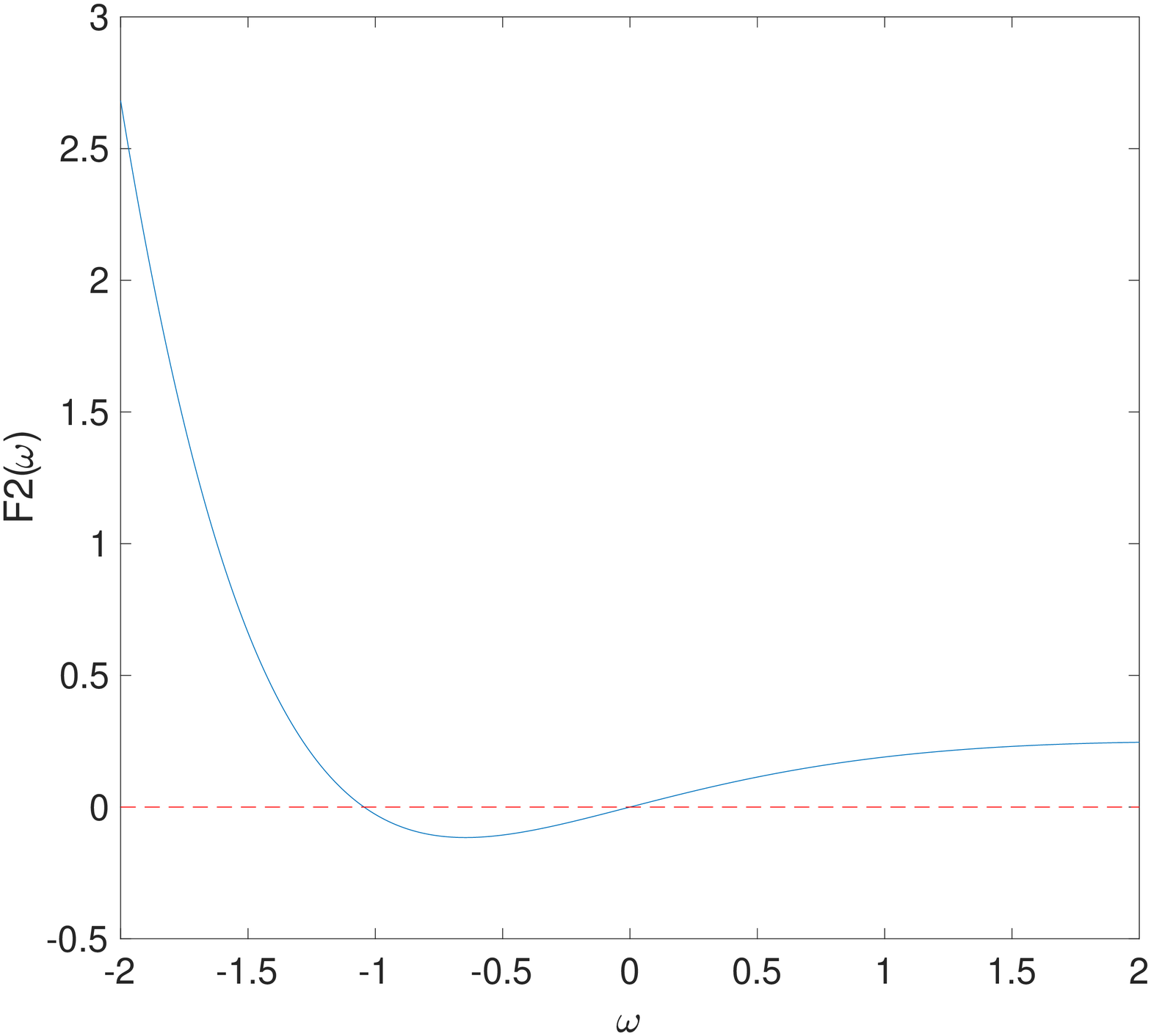}}
\caption{The plot of $F_1(u)$ and $F_1(u)$ . }
\label{fig:f1}
\end{figure}

We can proceed similarly for 
the system of equations \eqref{sys:III_b}. From its first equation   we see 
\begin{equation}
z=\frac{\om \frac{\widehat R_{1,N}^\theta}{\widehat R_{0,N}^\theta}-\th B(\om)}{\om-B(\om)}\,. \label{e.3.27}
\end{equation}
Substituting it to the second equation of \eqref{sys:III_b} we  have 
\begin{eqnarray}
&&\om 
\left(\frac{\widehat R_{1,N}^\theta}{\widehat R_{0,N}^\theta}
 -\theta\right)^2+\om \left(\om \frac{\widehat R_{1,N}^\theta}{\widehat R_{0,N}^\theta}
-\theta B(\om)\right)^2\nonumber\\
&&\qquad\qquad +B(\om) \left[\theta^2(\om-B(\om))^2-\left(\om\frac{\widehat R_{1,N}^\theta}{\widehat R_{0,N}^\theta}-\th B(\om)\right)^2\right]  -\om (\om-B(\om))^2
\frac{\widehat R_{2,N}^\theta}{\widehat R_{0,N}^\theta}=0\,.  \label{e.3.28} 
\end{eqnarray} 
 Denote the
left-hand side by $F_2(\om)$. Solving $F_2(\om)=0$ and
  substituting 
it into \eqref{e.3.32} yields  $\widehat  z_N$.
Notice that the 
quantities $\frac{\widehat R_{1,N}^\theta}{\widehat R_{0,N}^\theta}$ and $\frac{\widehat R_{2,N}^\theta}{\widehat R_{0,N}^\theta}$ appeared in \eqref{e.3.33} can  also be computed from real data. We simulate a sample of the process \eqref{e:MTOU} and plot the graphs $F_1(u)$ and $F_2(\om)$ in Figure \ref{fig:f1}.
 We take $\sigma=1$, $\alpha_1=0.1$, $\alpha_2=0.5$, $\beta_1=0.2$, $\beta_2=0.5$, $\theta=0.3$, $h=0.5$, $N=100,000$. It can be seen that since the case of  $u=0$ and $\om=0$ is excluded in Remark \ref{rmk:3.10}, there exists only one root for $F_1$ (or $F_2$).



\section{Numerical experiments}\label{sec:N}
To validate our estimation scheme discussed in Section \ref{sec:S},   we conduct some numerical experiments in this section. Table \ref{tab:T0Mean} and Table \ref{tab:T0Std} demonstrate the mean and standard deviation of the estimators $\widehat\alpha_{1,n,N}$ and $\widehat\alpha_{2,n,N}$ with $\sigma=1$  and  taking the order $n\in\{1,2,3,4,5,6,7\}$, $\theta=0$, $\beta_1=\beta_2=0$  through 1,000 sample paths. Here we set the simulation parameters as: $h=0.5$, $N=100,000$, $X_0=0$. Based on the numerical results, it can be seen that   the estimators have good consistency  and the estimators corresponding the order $n=2,3$ are recommended.
By using the built-in function ``fsolve'' in Matlab to solve the system in \eqref{e:L/L},  for given parameters $h=0.5$,  $\theta=0.1$, $\beta_1=\beta_2=0$, and $\sigma=0.6$, we estimate the parameters $\alpha_1$ and $\alpha_2$ in Table \ref{tab:TheMean} and  show the standard deviation in Table \ref{tab:TheStd}.


\begin{SCtable}
\caption{Mean of the estimators $\widehat\alpha$  through 1,000 sample paths. The true parameters are setting as: $\alpha_1=0.02$, $\alpha_2=0.05$.}\label{tab:T0Mean} 
    \begin{tabular}{lllllllllllllllllll}\topline
\rowcolor{gray!10}&  n&&&&&& \\ \cline{2-8}  
\rowcolor{gray!10} Mean&1&2&3&4&5&6&7\\ \rowmidlineHR
\rowcolor{gray!10}$ \alpha_1$&0.0198&0.0199&0.0199&0.0200&0.0200&0.0200&0.0201\\
\rowcolor{gray!10}$\alpha_2$&0.0497&0.0497&0.0495&0.0496&0.0497&0.0499&0.0497\\
\bottomlinec
    \end{tabular}
\end{SCtable}
\begin{SCtable}
\caption{Standard deviation of the estimators $\widehat\alpha$  through 1,000 sample paths. The true parameters are setting as: $\alpha_1=0.02$, $\alpha_2=0.05$.}\label{tab:T0Std}
  \begin{tabular}{lllllllllllllllll}\topline
\rowcolor{gray!10}&  n&&&&&& \\ \cline{2-8}  
\rowcolor{gray!10} Std&1&2&3&4&5&6&7\\  \rowmidlineHR
\rowcolor{gray!10}$ \alpha_1$&0.0012&0.0011&0.0011&0.0011&0.0012&0.0013&0.0015\\
\rowcolor{gray!10} $\alpha_2$&0.0027&0.0025&0.0023&0.0023&0.0025&0.0026&0.0030\\
\bottomlinec
    \end{tabular}
\end{SCtable}
\begin{SCtable}
\caption{Mean of the estimators $\widehat\alpha$ through 1,000 sample paths. The true parameters are setting as: $\alpha_1=0.1$, $\alpha_2=0.2$, $\theta=0.1$}\label{tab:TheMean}
 \begin{tabular}{lllllllllllllllllll}\topline
 \rowcolor{gray!10} & N($\times 10^4$)&&& \\ \cline{2-10}  
\rowcolor{gray!10} Mean&0.8          &1.2\quad\quad\quad  &1.6\quad\quad\quad  &2.0\quad\quad\  \\ \rowmidlineHR
\rowcolor{gray!10}$ \alpha_1$&0.0981&0.0979&0.0977&0.0974\\
\rowcolor{gray!10} $\alpha_2$&0.1917&0.1911&0.1910&0.1908\\
\bottomlinec
    \end{tabular}
\end{SCtable}

\begin{SCtable}
\caption{Standard deviation of the estimators $\widehat\alpha$  through 1,000 sample paths. The true parameters are setting as: $\alpha_1=0.1$, $\alpha_2=0.2$, $\theta=0.1$}\label{tab:TheStd}
    \begin{tabular}{lllllllllllllllllll}\topline
   \rowcolor{gray!10}        & N($\times 10^4$)&&& \\ \cline{2-8}  
 \rowcolor{gray!10} Std&0.8          &1.2\quad\quad\quad  &1.6\quad\quad\quad  &2.0\quad\quad\quad \ \\ \rowmidlineHR
 \rowcolor{gray!10}$ \alpha_1$&0.0094&0.0074&0.0065&0.0056\\
 \rowcolor{gray!10} $\alpha_2$&0.0150&0.0122&0.0103&0.0090\\
\bottomlinec
    \end{tabular}
\end{SCtable}
\section{Conclusion}\label{sec:C}
We conclude the paper here. In this paper, we have proposed the stationary moment estimators for the two-regime threshold OU process. Our approach can be extended to more threshold diffusion processes, including the threshold square-root process, where $X$ is a positive process almost surely with the diffusion term $\sigma(x)=\sum_{i=1}^m\sigma_i\sqrt x I(\theta_{i-1}<x\le\theta_i)$, $0=\theta_0<\theta_1<\theta_2<\cdots<\theta_m=\infty$. In the multi threshold OU  case, the stationary density is given by
\begin{equation*}
\psi(x)=\sum_{i=1}^mk_i\exp\left(\frac{-\alpha_ix^2+\beta_ix}{\sigma_i^2}\right)I(\theta_{i-1}<x\le \theta_i),
\end{equation*}
with $k_i$ determined by $\int_{-\infty}^\infty\psi(x)dx$ and $\sigma_i^2\psi(\theta_i-)=\sigma^2_{i+1}\psi(\theta_i+)$, $i=1,\ldots,m-1$. Notice that $\psi(x)$ may be not continuous at the point $\theta_i$.  
In addition, our estimation approach may be extended to estimate $\alpha_i$, $\beta_i$, $\theta$, and $\sigma$ simultaneously. For the related reference, we mention the recent work in \citet{cheng2020generalized}. They employed the ergodic theorem for $X_{t_k}-X_{t_{k-1}}$, and derived its characteristic function under the stationary distribution. For the threshold process, it is more difficult to conduct these because of the nonlinear term of the threshold process. The problem of estimating $\alpha_i$, $\beta_i$, $\theta$, and $\sigma$ simultaneously will be studied in a  future work.

\appendix
\section{Proof of Proposition \ref{Pro:spec}}\label{A:spectral}
We compute the transition probability by the spectral expansion method in \citet{linetsky2005transition}.  The proof is similar to Theorem 3.2 in \citet{decamps2006self} and  Proposition 3.1 in \citet{wang2015skew}.   So we just show the main computation  procedure here. For more details we refer the reader to Proposition 3.1 in \citet{wang2015skew}.
The spectral expansion of the density is written as
\begin{equation}
p_t(x,y)=m(y)\sum_{n=1}^{\infty}\exp(-\lambda t)\varphi_n(x)\varphi_n(y),\label{e.A.1} 
\end{equation}
where $\varphi_n(x)$ is the normalized eigenfunction associated to $\lambda_n$. It is well-known that 
\begin{equation*}
\xi(x,\lambda)=\exp\left( z_1^2/4\right)D_v(-z_1)
\end{equation*}
and 
\begin{equation*}
\eta(x,\lambda)=\exp\left( z_2^2/4\right)D_v(z_2)
\end{equation*}
are the solutions with continuous scale derivatives to the following Strum-Liouville equation
\begin{equation*}
\frac{1}{2}\sigma^2u^{\prime\prime}(x)+(\beta_1-\alpha_1x)u^\prime(x)=-\lambda u(x),\quad x\le\theta,
\end{equation*}
and 
\begin{equation*}
\frac{1}{2}\sigma^2u^{\prime\prime}(x)+(\beta_2-\alpha_2x)u^\prime(x)=-\lambda u(x),\quad x\ge\theta,
\end{equation*}
respectively. 

The Wronskian is given by
\begin{align*}
\omega(\lambda)&=\xi(\theta, \lambda)\frac{\eta^\prime(\theta, \lambda)}{s(\theta)}-\eta(\theta, \lambda)\frac{\xi^\prime(\theta, \lambda)}{s(\theta)},
\end{align*}
where $\eta^\prime(\theta,\lambda)=\frac{\partial\eta(x,\lambda)}{\partial x}\bigg|_{x=\theta}$ and $\xi^\prime(\theta,\lambda)=\frac{\partial\xi(x,\lambda)}{\partial x}\bigg|_{x=\theta}$. Noticing that the Hermite function $H_v(z)$ satisfies the recurrence relation \citep[see][Page 289]{MR0174795} as
\begin{equation*}
\frac{\partial H_v(z)}{\partial z}=2v  H_{v-1}(z),
\end{equation*}
 we get functions $\omega(\lambda)$ and $\varphi_n(x)$ in \eqref{e:w} and \eqref{e:varphi} respectively. Thus, the proof is completed .

\section{Computation of the asymptotic covariances}
In this section we compute the covariance $\sigma  $ 
 in Theorems \ref{asy:1},   \ref{t.3.7} and \ref{t.3.10} 
in  details by using the invariant measure $\psi_3$ given by
\eqref{e.general_psi} and  by the transition probability density function. We give a general formula. 
For any functions $f$ and $g$ we denote $\langle f\rangle =\int_\RR f(x)dx$ and $\langle f,g\rangle=\langle  fg\rangle $.
 Let $p_t(x,y)$ be the transition density of \eqref{e:MTOU}
 which is also  given by \eqref{e.A.1}. 
 Define $P_tf(x)=\int_\RR p_t(x, y) f(y) dy$.
 Then,  we  have   $\EE \left[ g(\widetilde X_{kh})| \widetilde X_0
 \right] =P_{kh}g(\widetilde X_0)$.  
For any two functions $f, g:\RR\rightarrow \RR$ if the following covariance is convergent, then it can be computed as 
\begin{align*}
\sigma(f,g)  &=\Cov\left [ f(\widetilde X_0)\,, g(\widetilde X_0)  \right]+ \sum_{k=1}^\infty \Cov\left[ f(\widetilde X_0)\,,  g( \widetilde X_{kh})\right]
+ \sum_{k=1}^\infty \Cov\left[g(\widetilde X_0)\,,   f( \widetilde X_{kh})\right]\\
&=\Cov\left [ f(\widetilde X_0)\,, g(\widetilde X_0)  \right]+ \sum_{k=1}^\infty \left\{\EE\left[ f(\widetilde X_0)  g( \widetilde X_{kh})\right]-\EE\left[ f(\widetilde X_0)  \right]\EE\left[g( \widetilde X_{kh})\right]\right\}\\
&\qquad\qquad + \sum_{k=1}^\infty \left\{\EE\left[g(\widetilde X_0)  f ( \widetilde X_{kh})\right]-\EE\left[g(\widetilde X_0)  \right]\EE\left[ f( \widetilde X_{kh})\right]\right\}
 \\
 &=\Cov\left [ f(\widetilde X_0)\,, g(\widetilde X_0)  \right]+ \sum_{k=1}^\infty \left\{\EE\left[ f(\widetilde X_0)  (P_{kh}g)( \widetilde X_{0})\right]-\EE\left[ f(\widetilde X_0)  \right]\EE\left[(P_{kh}g)( \widetilde X_{0})\right]\right\}\\
&\qquad\qquad + \sum_{k=1}^\infty \left\{\EE\left[g(\widetilde X_0)  (P_{kh}f )( \widetilde X_{0})\right]-\EE\left[g(\widetilde X_0)  \right]\EE\left[(P_{kh} f)( \widetilde X_{0})\right]\right\}\,. 
\end{align*}
Denote $\psi_f :=\sum_{k=1}^\infty  P_{kh} f $ and $\widehat \psi_f=
\psi_f-\langle \psi_3, \psi_f\rangle$.   Then, we have
\begin{equation}\label{e.B.1}
  \psi_f =(I-P_{h})^{-1}P_{h}  f \,, \quad
 \widehat  \psi_f=  (I-P_{h})^{-1}P_{h}  f- 
 \langle \psi_3, (I-P_{h})^{-1}P_{h}  f\rangle\,. 
\end{equation}
We also denote $\widehat  f=f-\langle \psi_3, f\rangle=\langle \psi_3  f \rangle$. 
With these notations, we have    
\begin{align}
\si(f,g) 
 &=\langle \psi_3   f g \rangle -\langle \psi_3  f \rangle \langle \psi_3 g\rangle  +  \langle \psi_3   g \psi_f\rangle  -
   \langle \psi_3   g  \rangle \langle \psi_3   \psi_f  \rangle  
+ \langle \psi_3    f \psi_g \rangle  -
   \langle \psi_3   f  \rangle \langle \psi_3   \psi_g  \rangle\nonumber\\
 &=\langle \psi_3  \widehat  f \widehat  g \rangle   +  \langle \psi_3   g \widehat \psi_f\rangle  
+ \langle   \psi_3    f \widehat  \psi_g \rangle  \,.\label{e.B.2}   
\end{align}
%
\noindent\textbf{Declarations of interest:} none.

\bibliography{\jobname}

\begin{thebibliography}{34}
\providecommand{\natexlab}[1]{#1}
\providecommand{\url}[1]{\texttt{#1}}
\expandafter\ifx\csname urlstyle\endcsname\relax
  \providecommand{\doi}[1]{doi: #1}\else
  \providecommand{\doi}{doi: \begingroup \urlstyle{rm}\Url}\fi

\bibitem[Bass and Pardoux(1987)]{bass1987uniqueness}
R.~F. Bass and \'{E}. Pardoux.
\newblock Uniqueness for diffusions with piecewise constant coefficients.
\newblock \emph{Probab. Theory Related Fields}, 76\penalty0 (4):\penalty0
  557--572, 1987.
\newblock ISSN 0178-8051.
\newblock \doi{10.1007/BF00960074}.
\newblock URL \url{https://doi.org/10.1007/BF00960074}.

\bibitem[Brockwell et~al.(1991)Brockwell, Hyndman, and
  Grunwald]{brockwell1991continuous}
Peter~J. Brockwell, Rob~J. Hyndman, and Gary~K. Grunwald.
\newblock Continuous time threshold autoregressive models.
\newblock \emph{Statist. Sinica}, 1\penalty0 (2):\penalty0 401--410, 1991.
\newblock ISSN 1017-0405.

\bibitem[Brockwell et~al.(2007)Brockwell, Davis, and
  Yang]{brockwell2007continuous}
Peter~J. Brockwell, Richard~A. Davis, and Yu~Yang.
\newblock Continuous-time {G}aussian autoregression.
\newblock \emph{Statistica Sinica}, 17\penalty0 (1):\penalty0 63--80, 2007.
\newblock ISSN 10170405, 19968507.
\newblock URL \url{http://www.jstor.org/stable/26432511}.

\bibitem[Brockwell and Hyndman(1992)]{brockwell1992continuous}
P.J. Brockwell and R.J. Hyndman.
\newblock On continuous-time threshold autoregression.
\newblock \emph{International Journal of Forecasting}, 8\penalty0 (2):\penalty0
  157 -- 173, 1992.
\newblock ISSN 0169-2070.
\newblock \doi{https://doi.org/10.1016/0169-2070(92)90116-Q}.
\newblock URL
  \url{http://www.sciencedirect.com/science/article/pii/016920709290116Q}.

\bibitem[Brooks et~al.(2011)Brooks, Gelman, Jones, and
  Meng]{brooks2011handbook}
Steve Brooks, Andrew Gelman, Galin~L. Jones, and Xiao-Li Meng, editors.
\newblock \emph{Handbook of {M}arkov chain {M}onte {C}arlo}.
\newblock Chapman \& Hall/CRC Handbooks of Modern Statistical Methods. CRC
  Press, Boca Raton, FL, 2011.
\newblock ISBN 978-1-4200-7941-8.
\newblock \doi{10.1201/b10905}.
\newblock URL \url{https://doi.org/10.1201/b10905}.

\bibitem[Browne and Whitt(1995)]{browne1995piecewise}
Sid Browne and Ward Whitt.
\newblock Piecewise-linear diffusion processes.
\newblock In \emph{Advances in queueing}, Probab. Stochastics Ser., pages
  463--480. CRC, Boca Raton, FL, 1995.

\bibitem[Buchholz(1969)]{MR0240343}
Herbert Buchholz.
\newblock \emph{The confluent hypergeometric function with special emphasis on
  its applications}.
\newblock Translated from the German by H. Lichtblau and K. Wetzel. Springer
  Tracts in Natural Philosophy, Vol. 15. Springer-Verlag New York Inc., New
  York, 1969.

\bibitem[Chan(1993)]{chan1993consistency}
K.~S. Chan.
\newblock Consistency and limiting distribution of the least squares estimator
  of a threshold autoregressive model.
\newblock \emph{Ann. Statist.}, 21\penalty0 (1):\penalty0 520--533, 1993.
\newblock ISSN 0090-5364.
\newblock \doi{10.1214/aos/1176349040}.
\newblock URL \url{https://doi.org/10.1214/aos/1176349040}.

\bibitem[Cheng et~al.(2020)Cheng, Hu, and Long]{cheng2020generalized}
Yiying Cheng, Yaozhong Hu, and Hongwei Long.
\newblock Generalized moment estimators for {$\alpha$}-stable
  {O}rnstein-{U}hlenbeck motions from discrete observations.
\newblock \emph{Stat. Inference Stoch. Process.}, 23\penalty0 (1):\penalty0
  53--81, 2020.
\newblock ISSN 1387-0874.
\newblock \doi{10.1007/s11203-019-09201-4}.
\newblock URL \url{https://doi.org/10.1007/s11203-019-09201-4}.

\bibitem[Chi et~al.(2017)Chi, Dong, and Wong]{chi2017option}
Zeyu Chi, Fangyuan Dong, and Hoi~Ying Wong.
\newblock Option pricing with threshold mean reversion.
\newblock \emph{Journal of Futures Markets}, 37\penalty0 (2):\penalty0
  107--131, 2017.
\newblock \doi{10.1002/fut.21795}.
\newblock URL \url{https://onlinelibrary.wiley.com/doi/abs/10.1002/fut.21795}.

\bibitem[Decamps et~al.(2006)Decamps, Goovaerts, and
  Schoutens]{decamps2006self}
Marc Decamps, Marc Goovaerts, and Wim Schoutens.
\newblock Self exciting threshold interest rates models.
\newblock \emph{Int. J. Theor. Appl. Finance}, 9\penalty0 (7):\penalty0
  1093--1122, 2006.
\newblock ISSN 0219-0249.
\newblock \doi{10.1142/S0219024906003937}.
\newblock URL \url{https://doi.org/10.1142/S0219024906003937}.

\bibitem[Ding et~al.(2020)Ding, Cui, and
  Wang]{doi:10.1080/14697688.2020.1781235}
Kailin Ding, Zhenyu Cui, and Yongjin Wang.
\newblock A markov chain approximation scheme for option pricing under skew
  diffusions.
\newblock \emph{Quantitative Finance}, 0\penalty0 (0):\penalty0 1--20, 2020.
\newblock \doi{10.1080/14697688.2020.1781235}.
\newblock URL \url{https://doi.org/10.1080/14697688.2020.1781235}.

\bibitem[Gairat and Shcherbakov(2017)]{gairat2017density}
Alexander Gairat and Vadim Shcherbakov.
\newblock Density of skew {B}rownian motion and its functionals with
  application in finance.
\newblock \emph{Math. Finance}, 27\penalty0 (4):\penalty0 1069--1088, 2017.
\newblock ISSN 0960-1627.
\newblock \doi{10.1111/mafi.12120}.
\newblock URL \url{https://doi.org/10.1111/mafi.12120}.

\bibitem[Hu and Song(2013)]{hu2013parameter}
Yaozhong Hu and Jian Song.
\newblock Parameter estimation for fractional {O}rnstein-{U}hlenbeck processes
  with discrete observations.
\newblock In \emph{Malliavin calculus and stochastic analysis}, volume~34 of
  \emph{Springer Proc. Math. Stat.}, pages 427--442. Springer, New York, 2013.
\newblock \doi{10.1007/978-1-4614-5906-4_19}.
\newblock URL \url{https://doi.org/10.1007/978-1-4614-5906-4_19}.

\bibitem[Hu et~al.(2015)Hu, Lee, Lee, and Song]{hu2015parameter}
Yaozhong Hu, Chihoon Lee, Myung~Hee Lee, and Jian Song.
\newblock Parameter estimation for reflected {O}rnstein-{U}hlenbeck processes
  with discrete observations.
\newblock \emph{Stat. Inference Stoch. Process.}, 18\penalty0 (3):\penalty0
  279--291, 2015.
\newblock ISSN 1387-0874.
\newblock \doi{10.1007/s11203-014-9112-7}.
\newblock URL \url{https://doi.org/10.1007/s11203-014-9112-7}.

\bibitem[Jiang et~al.(2018)Jiang, Song, and Wang]{jiang2018pricing}
Yiming Jiang, Shiyu Song, and Yongjin Wang.
\newblock Pricing {E}uropean vanilla options under a jump-to-default threshold
  diffusion model.
\newblock \emph{J. Comput. Appl. Math.}, 344:\penalty0 438--456, 2018.
\newblock ISSN 0377-0427.
\newblock \doi{10.1016/j.cam.2018.04.039}.
\newblock URL \url{https://doi.org/10.1016/j.cam.2018.04.039}.

\bibitem[Karlin and Taylor(1981)]{karlin1981second}
Samuel Karlin and Howard~M. Taylor.
\newblock \emph{A second course in stochastic processes}.
\newblock Academic Press, Inc. [Harcourt Brace Jovanovich, Publishers], New
  York-London, 1981.

\bibitem[Kutoyants(2012)]{kutoyants2012identification}
Yury~A. Kutoyants.
\newblock On identification of the threshold diffusion processes.
\newblock \emph{Ann. Inst. Statist. Math.}, 64\penalty0 (2):\penalty0 383--413,
  2012.
\newblock ISSN 0020-3157.
\newblock \doi{10.1007/s10463-010-0318-1}.
\newblock URL \url{https://doi.org/10.1007/s10463-010-0318-1}.

\bibitem[Lebedev(1965)]{MR0174795}
N.~N. Lebedev.
\newblock \emph{Special functions and their applications}.
\newblock Revised English edition. Translated and edited by Richard A.
  Silverman. Prentice-Hall, Inc., Englewood Cliffs, N.J., 1965.

\bibitem[Lejay and Pigato(2020)]{doi:10.1111/sjos.12417}
Antoine Lejay and Paolo Pigato.
\newblock Maximum likelihood drift estimation for a threshold diffusion.
\newblock \emph{Scandinavian Journal of Statistics}, 47\penalty0 (3):\penalty0
  609--637, 2020.
\newblock \doi{10.1111/sjos.12417}.
\newblock URL \url{https://onlinelibrary.wiley.com/doi/abs/10.1111/sjos.12417}.

\bibitem[Linetsky(2005)]{linetsky2005transition}
Vadim Linetsky.
\newblock On the transition densities for reflected diffusions.
\newblock \emph{Adv. in Appl. Probab.}, 37\penalty0 (2):\penalty0 435--460,
  2005.
\newblock ISSN 0001-8678.
\newblock \doi{10.1239/aap/1118858633}.
\newblock URL \url{https://doi.org/10.1239/aap/1118858633}.

\bibitem[Meyn and Tweedie(2009)]{meyn2012markov}
Sean Meyn and Richard~L. Tweedie.
\newblock \emph{Markov chains and stochastic stability}.
\newblock Cambridge University Press, Cambridge, second edition, 2009.
\newblock ISBN 978-0-521-73182-9.
\newblock \doi{10.1017/CBO9780511626630}.
\newblock URL \url{https://doi.org/10.1017/CBO9780511626630}.
\newblock With a prologue by Peter W. Glynn.

\bibitem[Meyn and Tweedie(1993)]{meyn1993stability}
Sean~P. Meyn and R.~L. Tweedie.
\newblock Stability of {M}arkovian processes. {II}. {C}ontinuous-time processes
  and sampled chains.
\newblock \emph{Adv. in Appl. Probab.}, 25\penalty0 (3):\penalty0 487--517,
  1993.
\newblock ISSN 0001-8678.
\newblock \doi{10.2307/1427521}.
\newblock URL \url{https://doi.org/10.2307/1427521}.

\bibitem[Siu(2016)]{siu2016self}
Tak~Kuen Siu.
\newblock A self-exciting threshold jump-diffusion model for option valuation.
\newblock \emph{Insurance Math. Econom.}, 69:\penalty0 168--193, 2016.
\newblock ISSN 0167-6687.
\newblock \doi{10.1016/j.insmatheco.2016.05.008}.
\newblock URL \url{https://doi.org/10.1016/j.insmatheco.2016.05.008}.

\bibitem[Siu et~al.(2006)Siu, Tong, and Yang]{siu2006option}
Tak~Kuen Siu, Howell Tong, and Hailiang Yang.
\newblock Option pricing under threshold autoregressive models by threshold
  {E}sscher transform.
\newblock \emph{J. Ind. Manag. Optim.}, 2\penalty0 (2):\penalty0 177--197,
  2006.
\newblock ISSN 1547-5816.
\newblock \doi{10.3934/jimo.2006.2.177}.
\newblock URL \url{https://doi.org/10.3934/jimo.2006.2.177}.

\bibitem[Stramer and Roberts(2007)]{stramer2007bayesian}
O.~Stramer and G.~O. Roberts.
\newblock On {B}ayesian analysis of nonlinear continuous-time autoregression
  models.
\newblock \emph{J. Time Ser. Anal.}, 28\penalty0 (5):\penalty0 744--762, 2007.
\newblock ISSN 0143-9782.
\newblock \doi{10.1111/j.1467-9892.2007.00549.x}.
\newblock URL \url{https://doi.org/10.1111/j.1467-9892.2007.00549.x}.

\bibitem[Stramer et~al.(1996)Stramer, Tweedie, and
  Brockwell]{stramer1996existence}
O.~Stramer, R.~L. Tweedie, and P.~J. Brockwell.
\newblock Existence and stability of continuous time threshold {ARMA}
  processes.
\newblock \emph{Statist. Sinica}, 6\penalty0 (3):\penalty0 715--732, 1996.
\newblock ISSN 1017-0405.

\bibitem[Su and Chan(2015)]{su2015quasi}
Fei Su and Kung-Sik Chan.
\newblock Quasi-likelihood estimation of a threshold diffusion process.
\newblock \emph{J. Econometrics}, 189\penalty0 (2):\penalty0 473--484, 2015.
\newblock ISSN 0304-4076.
\newblock \doi{10.1016/j.jeconom.2015.03.038}.
\newblock URL \url{https://doi.org/10.1016/j.jeconom.2015.03.038}.

\bibitem[Su and Chan(2017)]{su2017testing}
Fei Su and Kung-Sik Chan.
\newblock Testing for threshold diffusion.
\newblock \emph{J. Bus. Econom. Statist.}, 35\penalty0 (2):\penalty0 218--227,
  2017.
\newblock ISSN 0735-0015.
\newblock \doi{10.1080/07350015.2015.1073594}.
\newblock URL \url{https://doi.org/10.1080/07350015.2015.1073594}.

\bibitem[Tong(1983)]{MR717388}
Howell Tong.
\newblock \emph{Threshold models in nonlinear time series analysis}, volume~21
  of \emph{Lecture Notes in Statistics}.
\newblock Springer-Verlag, New York, 1983.
\newblock ISBN 0-387-90918-4.
\newblock \doi{10.1007/978-1-4684-7888-4}.
\newblock URL \url{https://doi.org/10.1007/978-1-4684-7888-4}.

\bibitem[van~der Vaart(1998)]{MR1652247}
A.~W. van~der Vaart.
\newblock \emph{Asymptotic statistics}, volume~3 of \emph{Cambridge Series in
  Statistical and Probabilistic Mathematics}.
\newblock Cambridge University Press, Cambridge, 1998.
\newblock ISBN 0-521-49603-9; 0-521-78450-6.
\newblock \doi{10.1017/CBO9780511802256}.
\newblock URL \url{https://doi.org/10.1017/CBO9780511802256}.

\bibitem[Wang et~al.(2015)Wang, Song, and Wang]{wang2015skew}
Suxin Wang, Shiyu Song, and Yongjin Wang.
\newblock Skew {O}rnstein-{U}hlenbeck processes and their financial
  applications.
\newblock \emph{J. Comput. Appl. Math.}, 273:\penalty0 363--382, 2015.
\newblock ISSN 0377-0427.
\newblock \doi{10.1016/j.cam.2014.06.023}.
\newblock URL \url{https://doi.org/10.1016/j.cam.2014.06.023}.

\bibitem[Zhuo and Menoukeu-Pamen(2017)]{zhuo2017efficient}
Xiaoyang Zhuo and Olivier Menoukeu-Pamen.
\newblock Efficient piecewise trees for the generalized skew {V}asicek model
  with discontinuous drift.
\newblock \emph{Int. J. Theor. Appl. Finance}, 20\penalty0 (4):\penalty0
  1750028, 34, 2017.
\newblock ISSN 0219-0249.
\newblock \doi{10.1142/S0219024917500285}.
\newblock URL \url{https://doi.org/10.1142/S0219024917500285}.

\bibitem[Zhuo et~al.(2017)Zhuo, Xu, and Zhang]{zhuo2017simple}
Xiaoyang Zhuo, Guangli Xu, and Haoyan Zhang.
\newblock A simple trinomial lattice approach for the skew-extended {CIR}
  models.
\newblock \emph{Math. Financ. Econ.}, 11\penalty0 (4):\penalty0 499--526, 2017.
\newblock ISSN 1862-9679.
\newblock \doi{10.1007/s11579-017-0192-1}.
\newblock URL \url{https://doi.org/10.1007/s11579-017-0192-1}.

\end{thebibliography}

\end{document}